\documentclass[a4paper, 12pt]{article}
\pdfoutput=1

\usepackage[T1]{fontenc}
\usepackage[utf8]{inputenc}

\usepackage{amsfonts,amsmath,amssymb,mathtools,dsfont,microtype}
\usepackage[dvipsnames]{xcolor}
\usepackage[noBBpl]{mathpazo}

\DeclareFontFamily{U} {cmr}{}

\DeclareFontShape{U}{cmr}{m}{n}{
    <-6> cmr5
    <6-7> cmr6
    <7-8> cmr7
    <8-9> cmr8
    <9-10> cmr9
    <10-12> cmr10
    <12-> cmr12}{}

\DeclareSymbolFont{Xcmr} {U} {cmr}{m}{n}
\DeclareMathSymbol{\Delta}{\mathord}{Xcmr}{'001}
\DeclareMathSymbol{\Upsilon}{\mathord}{Xcmr}{'007}
\DeclareMathSymbol{\Omega}{\mathord}{Xcmr}{'012}
\DeclareMathSymbol{\Gamma}{\mathord}{Xcmr}{'000}
\DeclareMathSymbol{\Theta}{\mathord}{Xcmr}{'002}
\DeclareMathSymbol{\Lambda}{\mathord}{Xcmr}{'003}
\DeclareMathSymbol{\Xi}{\mathord}{Xcmr}{'004}
\DeclareMathSymbol{\Pi}{\mathord}{Xcmr}{'005}
\DeclareMathSymbol{\Sigma}{\mathord}{Xcmr}{'006}
\DeclareMathSymbol{\Phi}{\mathord}{Xcmr}{'010}
\DeclareMathSymbol{\Psi}{\mathord}{Xcmr}{'011}

\usepackage[labelsep = period, labelfont = bf, justification = centering]{caption}
\usepackage{float,graphicx,subcaption}
\DeclareCaptionSubType*[roman]{figure}

\usepackage{booktabs,makecell}
\usepackage{enumitem,mdwlist}
\setlist[itemize]{topsep=0ex,itemsep=0ex,parsep=0.4ex}
\setlist[enumerate]{topsep=0ex,itemsep=0ex,parsep=0.4ex}

\usepackage{pgf,tikz,ifthen,calc}
\usepackage{tkz-euclide}
\tkzSetUpPoint[fill=black, size = 3pt]
\tkzSetUpLine[color=black, line width=0.6pt]
\tkzSetUpCircle[color=black, line width=0.6pt]

\usepackage[left = 2.5cm, right = 2.5cm, top = 2.5cm, bottom = 2.5cm, headsep = 12pt,
    headheight = 15pt]{geometry}
\usepackage{parskip}

\usepackage[hyphens]{url}
\usepackage[linktoc = all, hidelinks, colorlinks, unicode=true,
    pagebackref=true]{hyperref}
\usepackage[capitalise, compress, nameinlink, noabbrev]{cleveref}

\hypersetup{linkcolor={blue!70!black}, citecolor={green!70!black},
    urlcolor={blue!70!black}}

\usepackage{amsthm}
\theoremstyle{plain}
\newtheorem{thm}{Theorem}[section]
\newtheorem{claim}[thm]{Claim}
\newtheorem{cor}[thm]{Corollary}

\newtheorem{lem}[thm]{Lemma}

\newtheorem{construction}[thm]{Construction}
\newtheorem{prob}[thm]{Problem}
\newtheorem{fact}[thm]{Fact}

\theoremstyle{definition}
\newtheorem{defi}[thm]{Definition}
\newtheorem{rem}[thm]{Remark}

\crefname{subsection}{\S}{\S\S}
\newenvironment{proofclaim}[1][Proof of claim]{\begin{proof}[#1]}{\end{proof}}

\DeclareFontFamily{U}{matha}{\hyphenchar\font45}
\DeclareFontShape{U}{matha}{m}{n}{
    <5> <6> <7> <8> <9> <10> gen * matha
    <10.95> matha10 <12> <14.4> <17.28> <20.74> <24.88> matha12
}{}
\DeclareSymbolFont{matha}{U}{matha}{m}{n}

\DeclareMathSymbol{\specialuparrow}{\mathrel}{matha}{"D2}
\DeclareMathSymbol{\specialrightarrow}{\mathrel}{matha}{"D1}

\renewcommand*{\backref}[1]{}
\renewcommand*{\backrefalt}[4]{
    \ifcase #1 Not cited.\or \(\specialuparrow\)#2\else \(\specialuparrow\)#2\fi}

\renewcommand{\epsilon}{\varepsilon}
\renewcommand{\ge}{\geqslant}
\renewcommand{\le}{\leqslant}
\renewcommand{\geq}{\geqslant}
\renewcommand{\leq}{\leqslant}
\renewcommand{\emptyset}{\varnothing}
\renewcommand{\subset}{\subseteq}

\newcommand{\defn}[1]{\textcolor{Maroon}{\emph{#1}}}

\newcommand{\ex}{\text{ex}}

\title{Ordered Trees with Linear\\
    Off-Diagonal Ramsey Numbers}
\author{
    Lior Gishboliner\footnotemark[1]
    \and
    Artur Kravtsov\footnotemark[2]
    \and
    Xiangyu Li\footnotemark[3]
}
\date{}

\begin{document}
\maketitle

\renewcommand{\thefootnote}{\fnsymbol{footnote}}\urlstyle{same}
\footnotetext[1]{Department of Mathematics, University of Toronto, Canada. \emph{Email}:
    \href{mailto:lior.gishboliner@utoronto.ca}{\nolinkurl{lior.gishboliner@utoronto.ca}}.
    Research supported by an NSERC Discovery Grant.}
\footnotetext[2]{Department of Mathematics, University of Toronto, Canada. \emph{Email}:
    \href{mailto:artur.kravtsov@mail.utoronto.ca}{\nolinkurl{artur.kravtsov@mail.utoronto.ca}}.}
\footnotetext[3]{Department of Mathematics, University of Toronto, Canada. \emph{Email}:
    \href{mailto:xiangyuu.li@mail.utoronto.ca}{\nolinkurl{xiangyuu.li@mail.utoronto.ca}}.
    Research supported by a University of Toronto Excellence Award (UTEA).}
\renewcommand{\thefootnote}{\arabic{footnote}}

\begin{abstract}
    \noindent We characterize the following three classes of ordered trees:
    \begin{enumerate}
        \item trees \(T\) satisfying \(R_<(T,K_n)=O(n)\);
        \item \(\chi\)-bounding trees;
        \item Ramsey-good trees.
    \end{enumerate}
    In particular, we show that the first two classes coincide. The second
    characterization answers a question of Axenovich, Rollin and Ueckerdt, and the third
    resolves a problem of Balko and Poljak.

    In our proof, we also introduce several novel constructions of ordered graphs with
    sublinear independence number, giving lower bounds on \(R_<(F,K_n)\) for various
    ordered forests \(F\).
\end{abstract}

\section{Introduction}

An \defn{ordered graph} is a graph equipped with a linear order on its
vertices,\footnote{We picture the vertices as ordered from left to right and use terms
such as ``leftmost'' and ``rightmost''.} and copies are required to preserve this order.
We use the following notions from the extremal theory of ordered graphs.
\begin{itemize}
    \item For an ordered graph \(F\), the \defn{ordered extremal number}
    \(\operatorname{ex}_{<}(N,F)\) is the maximum number of edges in an \(F\)-free
    ordered graph on \(N\) vertices.

    \item An ordered forest \(F\) is \defn{\(\chi\)-bounding with constant \(C\)} if
    every \(F\)-free ordered graph has chromatic number at most \(C\). It is
    \defn{\(\chi\)-bounding} if such a constant exists.

    \item For ordered graphs \(G\) and \(H\), the \defn{ordered Ramsey number}
    \(R_{<}(G,H)\) is the least \(N\) such that every red-blue edge-coloring of the
    complete ordered graph on \(N\) vertices contains a red copy of \(G\) or a blue copy
    of \(H\).

    \item An ordered graph \(F\) on \(k\) vertices is called \defn{Ramsey-\(n\)-good} if
    \(R_<(F,K_n)= \linebreak (k-1)(n-1)+1\).\footnote{For any connected \(k\)-vertex
    graph \(F\), ordered or unordered, the lower bound always holds: the disjoint union
    of \(n-1\) cliques of size \(k-1\) contains neither a copy of \(F\) nor an
    independent set of size \(n\).} It is called \defn{Ramsey-good} if it is
    Ramsey-\(n\)-good for \(n \geq 2\).

    \item The \defn{interval chromatic number} \(\chi_<(H)\) of an ordered graph \(H\)
    is the minimum \(k\) such that \(V(H)\) can be split into \(k\) intervals, each of
    which is an independent set.
\end{itemize}

Our work primarily concerns ordered Ramsey numbers, whose systematic study was initiated
independently by Balko et al.~\cite{BCKK} and Conlon et al.~\cite{CFLS}. See
\cite{BJV,GJJ,GJS,GL,MubayiSuk} for recent work on the subject, and Balko's
survey~\cite{Balko_survey}.

\subsection{Linear Ramsey numbers and Ramsey-goodness}

In this paper, we study the off-diagonal Ramsey number \(R_<(F,K_n)\), where \(F\) is a
fixed ordered forest. Our starting point is the following problem.
\begin{prob}\label{prob:Ramsey-linear}
    Characterize the ordered forests \(F\) for which \(R_<(F,K_n)=O(n)\).
\end{prob}

This problem is closely related to the characterization of \(\chi\)-bounding ordered
forests. Indeed, if \(F\) is \(\chi\)-bounding with constant \(C\), then every
\(F\)-free ordered graph \(G\) has an independent set of size at least \(|V(G)|/C\), and
hence \(R_<(F,K_n)=O(n)\). The study of \(\chi\)-bounding ordered graphs was initiated
by Axenovich, Rollin and Ueckerdt~\cite{ARU}. A beautiful recent result of Bourneuf,
Cocquet, Tang and Thomass{\'e}~\cite{Interval_minor} shows that all ordered matchings
are \(\chi\)-bounding. A strengthening of this -- showing that ordered graphs forbidding
a fixed ordered matching as an induced subgraph are \(\chi\)-bounded (meaning that their
chromatic number is upper-bounded in terms of their clique number) -- was announced by
Briański, Davies and Walczak \cite{BDW}.

Our main result resolves both characterization problems for ordered trees and shows that
the two properties coincide. To state the characterization, we recall the following
terminology from~\cite{FKMV}.

\begin{defi}
    An ordered tree \(T\) is \defn{increasing} if it is an edge or is obtained from an
    increasing tree \(T'\) by adding either a rightmost leaf adjacent to the leftmost
    vertex of \(T'\), or a leftmost leaf adjacent to the rightmost vertex of \(T'\).
\end{defi}

Note that every increasing tree \(T\) has \(\chi_<(T) = 2\), and there is an edge
between the leftmost vertex of \(T\) and its rightmost vertex.

\begin{defi}
    A \defn{\(z\)-tree} is an ordered tree obtained from an increasing tree with
    leftmost vertex \(\ell\) and rightmost vertex \(r\) by adding two sets of leaves
    \(A,B\), with \(A<\ell\) and \(B>r\), and the edges \(ra\) for all \(a\in A\) and
    \(\ell b\) for all \(b\in B\).
\end{defi}

The sets \(A,B\) may be empty. If either set is empty, the resulting tree is increasing.
See Figure~\ref{fig:z-tree} for an illustration of these definitions.
\begin{figure}[ht]
    \centering
    \begin{tikzpicture}[
        vertex/.style={circle,fill=black,inner sep=2pt},
        line cap=round
        ]
        \begin{scope}
            \node[vertex] (t1) at (0,0) {};
            \node[vertex] (t2) at (0.7,0) {};
            \node[vertex] (t3) at (1.4,0) {};
            \node[vertex] (t4) at (2.1,0) {};

            \draw[thick] (t2) to[bend left=35] (t3);
            \draw[thick] (t1) to[bend left=40] (t3);
            \draw[thick] (t1) to[bend left=50] (t4);
        \end{scope}

        \begin{scope}[xshift=4cm]
            \node[vertex] (z1) at (0,0) {};
            \node[vertex] (z2) at (0.7,0) {};
            \node[vertex] (z3) at (1.4,0) {};
            \node[vertex] (z4) at (2.1,0) {};
            \node[vertex] (z5) at (2.8,0) {};
            \node[vertex] (z6) at (3.5,0) {};
            \node[vertex] (z7) at (4.2,0) {};
            \node[vertex] (z8) at (4.9,0) {};

            \draw[thick] (z4) to[bend left=35] (z5);
            \draw[thick] (z3) to[bend left=40] (z5);
            \draw[thick] (z3) to[bend left=50] (z6);

            \draw[thick] (z1) to[bend left=55] (z6);
            \draw[thick] (z2) to[bend left=50] (z6);
            \draw[thick] (z3) to[bend left=50] (z7);
            \draw[thick] (z3) to[bend left=55] (z8);
        \end{scope}
    \end{tikzpicture}
    \caption{An increasing tree (left) and a \(z\)-tree (right).} \label{fig:z-tree}
\end{figure}
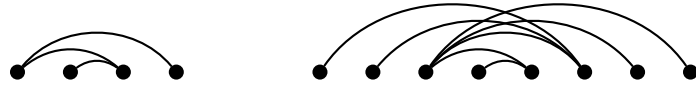
\begin{defi}[interval cut vertex]\label{def:interval cut vertex}
    Let \(F\) be a connected ordered graph. An \defn{interval cut vertex} of \(F\) is a
    vertex \(v\) that is neither the first nor the last vertex and for which no edge
    \(uw \in E(F)\) satisfies \(u < v < w\).
\end{defi}

If \(v\) is an interval cut vertex of \(F\), let \(I := \{u \in V(F) : u \leq v\}\) and
\(J := \{u \in V(F) : u \geq v\}\). Then every edge of \(F\) is contained in \(I\) or
\(J\), and both \(F[I]\) and \(F[J]\) are connected. Thus, by repeatedly splitting at
interval cut vertices, any connected ordered graph \(F\) can be uniquely decomposed into
intervals \(I_1 \leq \dots \leq I_m\) satisfying \(V(F) = I_1 \cup \dots \cup I_m\);
\(|I_i \cap I_{i+1}| = 1\) for \(1 \leq i < m\); \(I_i \cap I_j = \emptyset\) for \(j-i
\geq 2\); and \(F[I_i]\) has no interval cut vertices for every \(i \in [m]\). The
intervals \(I_1,\dots,I_m\) are called the \defn{segments} of \(F\).

Our main result can now be stated as follows.
\begin{thm}\label{thm:charac}
    For an ordered tree \(T\), the following are equivalent:
    \begin{enumerate}
        \item \(T\) is \(\chi\)-bounding.
        \item \(R_<(T,K_n) = O(n)\).
        \item For every segment \(I\) of \(T\), the induced subtree \(T[I]\) is a
        \(z\)-tree.
    \end{enumerate}
\end{thm}

This answers a question of Axenovich et al.~\cite{ARU}, who asked to characterize the
\(\chi\)-bounding ordered trees. The more general problem for ordered forests remains
open (see Problem~\ref{prob:Ramsey-linear}). We also note that for general ordered
forests, \(R_<(F,K_n)=O(n)\) need not imply that \(F\) is \(\chi\)-bounding; see
Section~\ref{subsec:additional} below.

The implication \(1 \Rightarrow 2\) in Theorem~\ref{thm:charac} is immediate. Every
\(z\)-tree \(T\) has \(\chi_<(T)=2\), and F\"uredi, Kostochka, Mubayi and
Verstra\"{e}te~\cite{FKMV} proved that if \(T\) has \(k \geq 2\) vertices, then
\(\ex_<(n,T)=(k-2)n-\binom{k-1}{2}\) for all \(n \geq k\). As we shall see, \(3
\Rightarrow 1\) follows easily from this bound. Thus, our main contribution is proving
\(2 \Rightarrow 3\).

A stronger property than satisfying \(R_<(F,K_n)=O(n)\) is Ramsey-goodness. In the
unordered setting, the notion of Ramsey-goodness goes back to Burr and
Erd\H{o}s~\cite{BE} and has been thoroughly studied since. A classical result of
Chv\'atal~\cite{Chvatal} states that every unordered tree is Ramsey-good. The study of
Ramsey-goodness for ordered graphs was initiated by Balko and Poljak~\cite{BP}. We prove
their conjecture by characterizing the Ramsey-good ordered trees.
\begin{defi}
    A \defn{one-sided star} is an ordered star whose center is the first or last vertex.
\end{defi}

\begin{thm}\label{thm:Ramsey-good}
    An ordered tree \(T\) is Ramsey-good if and only if, for every segment \(I\) of
    \(T\), the induced subgraph \(T[I]\) is a one-sided star.
\end{thm}

\subsection{Tight bounds for 4-vertex ordered trees}\label{subsec:4-vertex}
\begin{figure}[h]
    \centering
    \begin{tikzpicture}[
        vertex/.style={circle,fill=black,inner sep=2pt},
        every label/.style={label distance=3pt},
        line cap=round
        ]
        \begin{scope}
            \node[vertex,label=below:\(1\)] (c1) at (0,0) {};
            \node[vertex,label=below:\(2\)] (c2) at (0.7,0) {};
            \node[vertex,label=below:\(3\)] (c3) at (1.4,0) {};
            \node[vertex,label=below:\(4\)] (c4) at (2.1,0) {};

            \draw[thick] (c1) to[bend left=35] (c2);
            \draw[thick] (c3) to[bend left=35] (c4);
            \draw[thick] (c1) to[bend left=50] (c4);
        \end{scope}

        \begin{scope}[xshift=4cm]
            \node[vertex,label=below:\(1\)] (a1) at (0,0) {};
            \node[vertex,label=below:\(2\)] (a2) at (0.7,0) {};
            \node[vertex,label=below:\(3\)] (a3) at (1.4,0) {};
            \node[vertex,label=below:\(4\)] (a4) at (2.1,0) {};

            \draw[thick] (a1) to[bend left=35] (a2);
            \draw[thick] (a2) to[bend left=35] (a3);
            \draw[thick] (a1) to[bend left=50] (a4);
        \end{scope}

        \begin{scope}[xshift=8cm]
            \node[vertex,label=below:\(1\)] (b1) at (0,0) {};
            \node[vertex,label=below:\(2\)] (b2) at (0.7,0) {};
            \node[vertex,label=below:\(3\)] (b3) at (1.4,0) {};
            \node[vertex,label=below:\(4\)] (b4) at (2.1,0) {};

            \draw[thick] (b1) to[bend left=55] (b4);
            \draw[thick] (b1) to[bend left=38] (b3);
            \draw[thick] (b2) to[bend left=38] (b4);
        \end{scope}

        \begin{scope}[xshift=12cm]
            \node[vertex,label=below:\(1\)] (d1) at (0,0) {};
            \node[vertex,label=below:\(2\)] (d2) at (0.7,0) {};
            \node[vertex,label=below:\(3\)] (d3) at (1.4,0) {};
            \node[vertex,label=below:\(4\)] (d4) at (2.1,0) {};

            \draw[thick] (d1) to[bend left=35] (d2);
            \draw[thick] (d1) to[bend left=60] (d3);
            \draw[thick] (d2) to[bend left=60] (d4);
        \end{scope}
    \end{tikzpicture}
    \caption{The ordered paths \(F_1,F_2,F_3,F_4\), from left to right.}
    \label{fig:prime-paths}
\end{figure}
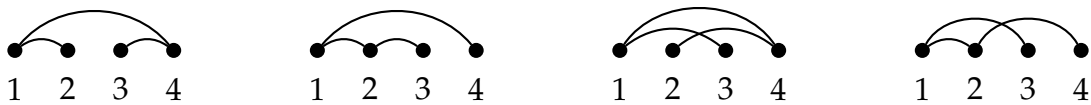

As noted above, the main difficulty in proving Theorem~\ref{thm:charac} is implication
\(2 \Rightarrow 3\). The proof proceeds by establishing superlinear lower bounds
\(R_<(T,K_n)=\omega(n)\) for several small ordered trees \(T\), and then showing that an
ordered tree with no interval cut vertices is a \(z\)-tree if it avoids all of these
trees. As a byproduct, we obtain tight bounds on \(R_<(T,K_n)\) for every four-vertex
ordered tree \(T\). Consider the four ordered trees \(F_1,F_2,F_3,F_4\) in
Figure~\ref{fig:prime-paths}.

\begin{thm}\label{thm:F1,F2}
    For the ordered trees \(F_1 = \{12,34,14\}\) and \(F_2 = \{12,23,14\}\), it holds
    that
    \[
        R_<(F_1,K_n) = \Theta(n \log n)
        \quad \text{and} \quad
        R_<(F_2,K_n) = \Theta(n \log n).
    \]
\end{thm}

We state these bounds together because both \(F_1\) and \(F_2\) are avoided by the same
construction, with a different vertex ordering for each (see
Construction~\ref{const:F1,F2}).

\begin{thm}\label{thm:F3}
    For the ordered tree \(F_3 = \{13,24,14\}\), it holds that
    \[
        R_<(F_3,K_n)
        =\Theta\left(\frac{n\log n}{\log\log n}\right).
    \]
\end{thm}

The lower bound in Theorem~\ref{thm:F3} is proved using a \defn{twisted hypercube} built
from random matchings (see Construction~\ref{const:F3}).

\begin{thm}\label{thm:F4}
    For the ordered tree \(F_4 = \{12,13,24\}\), it holds that
    \[
        R_<(F_4,K_n) = \Theta(n\log n).
    \]
\end{thm}

Together, Theorems~\ref{thm:charac}, \ref{thm:F1,F2}, \ref{thm:F3} and \ref{thm:F4}
determine \(R_<(F,K_n)\) up to a constant factor for every four-vertex ordered forest
\(F\). Indeed, for an ordered graph \(F\), let \(\overleftarrow{F}\) be the ordered
graph obtained from \(F\) by reversing the vertex order. Clearly, \(R_<(F,K_n) =
R_<(\overleftarrow{F},K_n)\). Thus, if \(F\) is one of
\(F_1,F_2,\overleftarrow{F_2},F_3,F_4,\overleftarrow{F_4}\), then \(R_<(F,K_n)\) is
determined by Theorems~\ref{thm:F1,F2}, \ref{thm:F3} and \ref{thm:F4}. In all other
cases, \(R_<(F,K_n) = O(n)\) (and thus \(R_<(F, K_n) = \Theta(n)\) if \(F\) has edges)
by Theorem~\ref{thm:charac}.

We end this section with one more small ordered tree, this time on 5 vertices, which
will be needed in the proof of Theorem~\ref{thm:charac}. Let \(R =\{13,23,14,25\}\) (see
the right graph in Figure \ref{fig:Q,R}).
\begin{thm}\label{thm:R}
    For the ordered tree \(R =\{13,23,14,25\}\), it holds that
    \[
        \Omega\left(\frac{n\log n}{(\log\log n)^2}\right) \leq R_<(R,K_n)
        \leq O\left(
        \frac{n\log n}{\log \log n} \right).
    \]
\end{thm}

The construction used for the lower bound in Theorem~\ref{thm:R} (see
Construction~\ref{const:randomized-Tutte}) is inspired by the construction of
\(K_3\)-free graphs with high chromatic number due to Tutte (alias Blanche Descartes).
This construction was used by Axenovich et al.~\cite{ARU} in their study of
\(\chi\)-bounding in ordered graphs.

\subsection{Lower bounds via Davenport-Schinzel sequences}
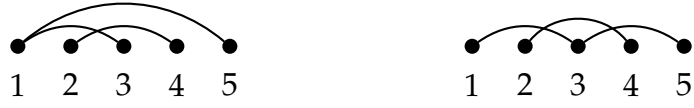
\begin{figure}
    \centering
    \begin{tikzpicture}[
        vertex/.style={circle,fill=black,inner sep=2pt},
        every label/.style={label distance=3pt},
        line cap=round
        ]

        \begin{scope}
            \node[vertex,label=below:\(1\)] (c1) at (0,0) {};
            \node[vertex,label=below:\(2\)] (c2) at (0.7,0) {};
            \node[vertex,label=below:\(3\)] (c3) at (1.4,0) {};
            \node[vertex,label=below:\(4\)] (c4) at (2.1,0) {};
            \node[vertex,label=below:\(5\)] (c5) at (2.8,0) {};

            \draw[thick] (c1) to[bend left=35] (c3);
            \draw[thick] (c2) to[bend left=35] (c4);
            \draw[thick] (c1) to[bend left=40] (c5);
        \end{scope}

        \begin{scope}[xshift=6cm]
            \node[vertex,label=below:\(1\)] (c1) at (0,0) {};
            \node[vertex,label=below:\(2\)] (c2) at (0.7,0) {};
            \node[vertex,label=below:\(3\)] (c3) at (1.4,0) {};
            \node[vertex,label=below:\(4\)] (c4) at (2.1,0) {};
            \node[vertex,label=below:\(5\)] (c5) at (2.8,0) {};

            \draw[thick] (c1) to[bend left=35] (c3);
            \draw[thick] (c2) to[bend left=50] (c4);
            \draw[thick] (c3) to[bend left=35] (c5);
        \end{scope}
    \end{tikzpicture}
    \caption{The two ordered forests in Theorems \ref{thm:F_0} and \ref{thm:13,35,24}.
    The left forest is \(F_0\).} \label{fig:F_0}
\end{figure}

This section concerns a superlinear lower bound on \(R_<(F,K_n)\) for two more ordered
forests. One of these is \(F_0 = \{13,15,24\}\); see the left graph in Figure
\ref{fig:F_0}. The fact that \(R_<(F_0,K_n) = \omega(n)\) will be used in the proof of
Theorem \ref{thm:charac}. We note that \(F_0\) is well-studied in the context of ordered
Tur\'an problems, due to its close relation to Davenport-Schinzel sequences.
Furthermore, the connection between Davenport-Schinzel sequences and ordered extremal
problems is well-known, going back to the work of F\"uredi and Hajnal
\cite{FurediHajnal}. Here we exploit this connection to show that \(R_<(F_0,K_n) =
\omega(n)\). Let \(\alpha(n)\) denote the inverse Ackermann function. This function
tends to infinity with \(n\) \nolinebreak (but \nolinebreak very \nolinebreak slowly).

\begin{thm}\label{thm:F_0}
    For the ordered forest \(F_0 = \{13,15,24\}\), it holds that
    \[
        \Omega\big( n \sqrt{\alpha(n)} \big) \leq R_<(F_0,K_n) \leq
        O(n \alpha(n)).
    \]
\end{thm}

The upper bound in Theorem \ref{thm:F_0} follows easily from the known fact that every
\(n\)-vertex \(F_0\)-free ordered graph has at most \(O(n \alpha(n))\) edges; see
\cite{PettieTardosWalczak}. Indeed, this bound implies that every \(n\)-vertex
\(F_0\)-free ordered graph is \(O(\alpha(n))\)-degenerate and hence has chromatic number
\(O(\alpha(n))\). We also note that \(F_0\) was very recently \cite{PettieTardosWalczak}
shown to not be \(\chi\)-bounding, and Theorem \ref{thm:F_0} strengthens this result.

Using Davenport-Schinzel sequences, we also show that the ordered forest \(F :=
\{13,35,24\}\) has a superlinear Ramsey number (see the right graph in Figure
\ref{fig:F_0}). We note that this \(F\) is one of the minimal ordered forests whose
\(\chi\)-bounding status was left open in \cite{ARU}. Our superlinear lower bound on
\(R_<(F,K_n)\) shows that \(F\) is not \(\chi\)-bounding, resolving this minimal case.
The fact that \(F\) is not \(\chi\)-bounding answers a question of Aboulker, Aubian and
Charbit \cite{AAC}.
\begin{thm}\label{thm:13,35,24}
    For the ordered forest \(F = \{13,35,24\}\), it holds that
    \[
        R_<(F,K_n) \geq \Omega\big( n \sqrt{\alpha(n)} \big).
    \]
\end{thm}

\subsection{Additional positive results}\label{subsec:additional}
Our last set of results consists of three reductions for showing that certain ordered
forests have a linear Ramsey number or are \(\chi\)-bounding. We believe that these are
useful in the further study of Problem \ref{prob:Ramsey-linear} and its
\(\chi\)-bounding analogue. Moreover, we obtain examples of ordered forests \(F\) which
are Ramsey-linear (i.e., satisfy \(R_<(F,K_n) = O(n)\)) despite not being
\(\chi\)-bounding.

Our first reduction allows us to place forests  \(H_1,\dots,H_s\) between the leaves of
a one-sided star; if \(H_1,\dots,H_s\) are Ramsey-linear then so is the resulting
configuration.

\begin{thm}\label{thm:insertion-into-stars}
    Let \(s \geq 1\) and let \(H_1,\dots,H_s\) be ordered graphs, possibly empty. Let
    \(F\) be the ordered graph on
    \[
        x < V(H_1) < y_1 < V(H_2) < y_2 < \dots < V(H_s) < y_s
    \]
    with the edges \(xy_1,\dots,xy_s\) and all the edges of \(H_1,\dots,H_s\). If
    \(R_<(H_i,K_n) = O(n)\) for every \(i \in [s]\), then \(R_<(F,K_n) = O(n)\).
\end{thm}

In particular, Theorem \ref{thm:insertion-into-stars} shows that \(F = \{12,15,34\}\)
satisfies \(R_<(F,K_n) = O(n)\) (by taking \(s=2\) and \(H_2\) being an edge). Axenovich
et al.~\cite{ARU} showed that this \(F\) is not \(\chi\)-bounding. Indeed, they observed
that the natural ordering of shift graphs is \(F\)-free. Shift graphs are known to have
unbounded chromatic number.

Our second reduction is concerned with the following definition: For two ordered graphs
\(A,B\) with vertex-sets \(V(A) = \{a_1 < \dots < a_k\}\) and \(V(B) = \{b_1 < \dots <
b_{\ell}\}\), where \(k,\ell \geq 2\), let \(A \diamond B\) denote the ordered graph
consisting of disjoint copies of \(A,B\) with the vertex order \(a_1 < \dots < a_{k-1} <
b_1 < a_k < b_2 < \dots < b_{\ell}\). In other words, \(A \diamond B\) is obtained by
first taking disjoint copies of \(A,B\) where \(B\) is completely to the right of \(A\),
and then swapping the places of the last vertex of \(A\) and the first vertex of \(B\).

\begin{thm}\label{thm:one-vertex interlace}
    For every two ordered graphs \(A,B\) with \(|V(A)|,|V(B)| \geq 2\), it holds that
    \[
        R_<(A \diamond B, K_n) \leq 2R_<(A,K_n) + 2R_<(B,K_n)-3.
    \]
\end{thm}

In particular, if \(R_<(A,K_n),R_<(B,K_n) = O(n)\), then the same holds for \(A \diamond
B\). This implies, e.g., that the ordered forest \(F = \{12,14,35\}\) has \(R_<(F,K_n) =
O(n)\) (here we take \(A\) to be the star \(\{a_1a_2,a_1a_3\}\) with \(a_1<a_2<a_3\),
and \(B\) to be an edge).

Our last reduction gives an operation preserving the property of being
\(\chi\)-bounding. For an ordered graph \(F\) with vertices \(v_1 < \dots < v_k\), let
\(F^{\star}\) be obtained from \(F\) by adding vertices \(a,b,c\) with \(a < v_1 < \dots
< v_k < b < c\) and edges \(ab,v_kc\).
\begin{thm}\label{thm:chi-bounding reduction}
    Let \(F\) be an ordered graph with vertices \(v_1 < \dots < v_k\) such that \(v_1v_k
    \in E(F)\). Suppose that \(F\) is \(\chi\)-bounding with constant \(C\). Then
    \(F^{\star}\) is \(\chi\)-bounding with constant \(5C+4\).
\end{thm}
Note that if \(F\) is an edge then \(F^{\star} = \{14,23,35\}\). Hence, this ordered
forest is \(\chi\)-bounding. This is one of the minimal graphs whose \(\chi\)-bounding
status was left open in \cite{ARU}.

\paragraph{Paper Organization.} Section \ref{sec:prelim} contains some preliminary
lemmas. In Section \ref{sec:4-vertex} we prove all theorems appearing in Section
\ref{subsec:4-vertex}. In Section \ref{sec:DS} we use Davenport-Schinzel sequences to
prove Theorems \ref{thm:F_0} and \ref{thm:13,35,24}. The proofs of Theorems
\ref{thm:charac} and \ref{thm:Ramsey-good} appear in Sections \ref{sec:charac} and
\ref{sec:Ramsey-good}, respectively. Section \ref{sec:reductions} gives the proofs of
Theorems \ref{thm:insertion-into-stars}, \ref{thm:one-vertex interlace} and
\ref{thm:chi-bounding reduction}. Finally, Section \ref{sec:concluding remarks}
contains some additional remarks.

\paragraph{Related work.} After completing this paper, we learned of related unpublished
results of Bria{\'n}ski, Davies, Galewski, and Walczak. They independently obtained the
equivalence \(1 \Leftrightarrow 3\) in Theorem~\ref{thm:charac} and showed that
\(\{14,23,35\}\) is \(\chi\)-bounding whereas \(\{13,35,24\}\) is not. Furthermore, they
announced that \(\{12,14,35\}\) is not \(\chi\)-bounding. Together with
Theorem~\ref{thm:one-vertex interlace}, this gives another example of Ramsey-linearity
without \(\chi\)-boundedness. We thank them for sharing their results with us.

\section{Preliminary lemmas}\label{sec:prelim}
In this section we collect several auxiliary lemmas on ordered Ramsey numbers and
\(\chi\)-bounding. The following lemma allows us to restrict our attention to ordered
trees having no interval cut vertices.
\begin{lem}\label{lem:segment}
    Let \(F\) be an ordered graph, let \(v \in V(F)\), and set
    \[
        I := \{u \in V(F) : u \leq v\}
        \quad \text{and} \quad
        J := \{u \in V(F) : u \geq v\}.
    \]
    Suppose that \(F\) has no edges between \(I \setminus \{v\}\) and \(J \setminus
    \{v\}\), and put \(F_1 := F[I]\) and \(F_2 := F[J]\). Then:
    \begin{enumerate}
        \item If \(F_i\) is \(\chi\)-bounding with constant \(C_i\) for \(i=1,2\), then
        \(F\) is \(\chi\)-bounding with constant \(C_1 + C_2\).
        \item \(R_<(F, K_n) = R_<(F_1, K_n) + R_<(F_2, K_n) - 1\).
    \end{enumerate}
\end{lem}
\begin{proof}
    \hfill
    \begin{enumerate}
        \item Let \(G\) be an \(F\)-free ordered graph. Let \(S_1\) be the set of
        vertices \(x \in V(G)\) such that there is no copy of \(F_1\) in \(G\) in which
        \(x\) plays the role of \(v\). Similarly, let \(S_2\) be the set of \(x \in
        V(G)\) such that there is no copy of \(F_2\) in \(G\) in which \(x\) plays the
        role of \(v\). Then \(S_1 \cup S_2 = V(G)\), because otherwise \(G\) would
        contain a copy of \(F\). Also, \(G[S_i]\) is \(F_i\)-free, and so \(\chi(G[S_i])
        \leq C_i\) by the assumptions of the lemma. The claim follows.
        \item For the lower bound, take an ordered graph \(G_i\) on \(R_<(F_i,K_n) - 1\)
        vertices with no copy of \(F_i\) and with \(\alpha(G_i) < n\). Let \(G\) be
        obtained by adding a complete bipartite graph between \(G_1,G_2\) and letting
        \(V(G_1) < V(G_2)\) in the vertex order. Then \(\alpha(G) < n\) and \(G\) is
        \(F\)-free. Indeed, suppose by contradiction that \(G\) contains a copy of
        \(F\), and let \(x \in V(G)\) be the vertex playing the role of \(v\). If \(x
        \in V(G_1)\) then \(G_1\) contains a copy of \(F_1\), and if \(x \in V(G_2)\)
        then \(G_2\) contains a copy of \(F_2\). In either case we get a contradiction.

        For the upper bound, let \(G\) be an ordered graph on \(R_<(F_1, K_n) + R_<(F_2,
        K_n) - 1\) vertices with \(\alpha(G) < n\). Let \(S_1,S_2\) be defined as in
        Item 1, i.e., \(S_i\) is the set of vertices \(x \in V(G)\) such that there
        exists no copy of \(F_i\) where \(x\) plays the role of \(v\). Then \(|S_i| <
        R_<(F_i,K_n)\), since \(S_i\) contains neither a copy of \(F_i\) nor an
        independent set of size \(n\). Hence, \(|S_1| + |S_2| \leq R_<(F_1,K_n) +
        R_<(F_2,K_n) - 2 = v(G) - 1 < v(G)\), meaning that there exists \(x \in V(G)
        \setminus (S_1 \cup S_2)\). This gives a copy of \(F\) in \(G\), as required.
        \qedhere
    \end{enumerate}
\end{proof}

Lemma \ref{lem:segment} immediately implies the following corollaries:
\begin{cor}\label{cor:segment chi-bounding}

    \hfill
    \begin{enumerate}
        \item An ordered tree \(F\) is \(\chi\)-bounding if and only if \(F[I]\) is
        \(\chi\)-bounding for every segment \(I\) of \(F\).
        \item An ordered tree \(T\) is Ramsey-good if and only if \(T[I]\) is
        Ramsey-good for every segment \(I\) of \(T\).
    \end{enumerate}
\end{cor}
\begin{proof}
    Item 1 is immediate from Lemma \ref{lem:segment}. To prove Item 2, it suffices to
    show that when \(T\) has exactly two segments \(I \leq J\), the tree \(T\) is
    Ramsey-good if and only if \(T_1 := T[I]\) and \(T_2 := T[J]\) are Ramsey-good. Note
    that \(T_1\) and \(T_2\) are trees. First, if \(T_1,T_2\) are Ramsey-good then by
    Item 2 of Lemma \ref{lem:segment}, for all \(n \geq 2\) we have
    \[
        \begin{aligned}
            R_<(T,K_n)
            &= R_<(T_1,K_n) + R_<(T_2,K_n) - 1 \\
            &= (|T_1|-1)(n-1) + (|T_2|-1)(n-1) + 1\\
            &
            = (|T|-1)(n-1) + 1,
        \end{aligned}
    \]
    as \(|T| = |T_1| + |T_2| - 1\). Conversely, if \(T_i\) is not Ramsey-good then
    \(R_<(T_i,K_n) > (|T_i| - \nolinebreak 1)(n-1) + 1\) for some \(n\) (because every
    tree \(T'\) has \(R_<(T',K_n) \geq \linebreak (|T'|-1)(n- \nolinebreak 1)+1\) for
    all \(n\)). Now, the same calculation (via Lemma \ref{lem:segment}) shows that
    \(R_<(T,K_n) > (|T|-1)(n-1)+1\), and so \(T\) is not Ramsey-good.
\end{proof}

The next lemma gives an upper bound on \(R_<(F,K_n)\) in the case that \(F\) has a
partition into intervals \(I < J\) with exactly one edge going between \(I,J\).
\begin{lem}\label{lem:e(I,J) = 1}
    Let \(F\) be an ordered graph with \(V(F)=I\sqcup J\), where \(I<J\). Put
    \(F_1=F[I]\) and \(F_2=F[J]\), and suppose that \(e_F(I,J)=1\). Put \(g(n) :=
    R_<(F_1, K_n) + R_<(F_2, K_n)\). Then for all \(n_1,n_2 \geq 1\) with \(n = n_1 +
    n_2\), we have
    \begin{equation}
        \label{eq:e(I,J) = 1}
        R_{<}(F,K_n)
        \le
        g(n) +
        R_<(F, K_{n_1}) + R_<(F, K_{n_2}).
    \end{equation}
    Consequently, for \(k \geq 0\), if \(g(n) = O(n \log^k n)\) then \(R_<(F,K_n) = O(n
    \log^{k+1}n)\).
\end{lem}

Note that Lemma \ref{lem:e(I,J) = 1} in particular applies if the first or last vertex
of \(F\) is a leaf.
\begin{proof}
    The second part of the lemma follows from the first one by the Master Theorem (see
    \cite[Section 4.5]{clrs}). We prove the first part. Put
    \[
        N := R_<(F_1, K_n) + R_<(F_2, K_n) +
        R_<(F, K_{n_1}) + R_<(F, K_{n_2})
    \]
    (this is the RHS in \eqref{eq:e(I,J) = 1}), and suppose by contradiction that there
    exists an \(N\)-vertex \(F\)-free ordered graph \(G\) with \(\alpha(G) < n\).
    Partition \(V(G)\) into intervals \(X, Y\) with \(X<Y\) of sizes
    \[
        |X| = R_<(F_1, K_n) + R_<(F, K_{n_1})
        \quad \text{and} \quad
        |Y| = R_<(F_2, K_n) + R_<(F, K_{n_2}).
    \]
    (This is possible by the choice of \(N\).) Let \(ab\) be the unique edge between
    \(I\) and \(J\), with \(a\in I\) and \(b\in J\). Let \(X'\) be the set of all \(x
    \in X\) such that there exists a copy of \(F_1\) in \(G[X]\) in which \(x\) plays
    the role of \(a\). Similarly, let \(Y'\) be the set of all \(y \in Y\) such that
    there exists a copy of \(F_2\) in \(G[Y]\) in which \(y\) plays the role of \(b\).
    If \(E(X', Y') \neq \varnothing\), then there exists a copy of \(F\), contradicting
    our assumption. Thus, \(E(X', Y') = \varnothing\).

    Next, observe that \(G[X \setminus X']\) is \(F_1\)-free (by the definition of
    \(X')\). Since \(G\) has no independent set of size \(n\), it follows that \(\vert X
    \setminus X' \vert < R_<(F_1, K_n)\). By the same argument, \(\vert Y \setminus Y'
    \vert < R_<(F_2, K_n)\). By our choice of the sizes of \(X,Y\), we conclude that
    \(\vert X'\vert \geq R_<(F,K_{n_1})\) and \(\vert Y'\vert \ge R_<(F,K_{n_2})\).
    Thus, \(X'\) contains an independent set of size \(n_1\), and \(Y'\) contains an
    independent set of size \(n_2\). But \(E(X',Y') = \varnothing\), giving an
    independent set of size \(n_1 + n_2 = n\) in \(G\), a contradiction.
\end{proof}

\begin{rem}
    One can prove a variant of the second part of Lemma \ref{lem:e(I,J) = 1} without any
    assumptions on \(g(n)\), instead assuming that \(F\) is connected. Indeed, \(F\)
    being connected implies that \(F_1,F_2\) are also connected. For a connected graph
    \(H\) with at least two vertices, the Ramsey numbers \(R_<(H,K_n)\) have the
    regularity property \(R_<(H,K_{n_1}) + R_<(H,K_{n_2}) \leq R_<(H,K_{n_1+n_2})\). (To
    see this, take \(H\)-free graphs \(G_1,G_2\) with \(|V(G_i)| = R_<(H,K_{n_i}) - 1\)
    and \(\alpha(G_i) \leq n_i-1\), and take the disjoint union of \(G_1,G_2\) and
    another isolated vertex.) It follows that if \(F\) is connected and \(|I|,|J| \geq
    2\), then \(g(n_1) + \nolinebreak g(n_2) \leq g(n_1+n_2)\). This regularity of \(g\)
    allows one to derive from \eqref{eq:e(I,J) = 1} that \(R_<(F,K_n) = O(g(n) \cdot
    \log n + n)\). The case that \(|I|=1\) or \(|J|=1\) can be handled similarly to
    obtain the same conclusion.
\end{rem}

\section{Constructions for 4-vertex ordered trees}\label{sec:4-vertex}

\subsection{Proof of Theorem \ref{thm:F1,F2}}

Recall that \(F_1 = \{12, 34, 14\}\) and \(F_2 = \{12, 23, 14\}\). To prove the lower
bounds in Theorem \ref{thm:F1,F2}, we will use the following construction.
\begin{construction}\label{const:F1,F2}
    Let \(G_k\) be the ordered graph whose vertices are the triples \((a,i,b) \in
    [k]^3\) with \(a < i < b\). The edges are defined as follows:

    \paragraph{Adjacency rule:} Two vertices \((a,i,b)\) and \((c,j,d)\) with \(i<j\)
    are adjacent if and only if \(i\le c<b\le j\).\footnote{In other words, \(u =
    (a,i,b)\) and \(v = (c,j,d)\) are adjacent if the points \((i,b)\) (i.e., the last
    two coordinates of \(u\)) interlace the points \((c,j)\) (i.e., the first two
    coordinates of \(v\)). Alternatively, one can also view a vertex \((a,i,b)\) as an
    interval \([a,b]\) with a distinguished interior point \(i\). Two vertices
    \((a,i,b)\) and \((c,j,d)\) with \(i<j\) are adjacent when the corresponding
    intervals \([a,b],[c,d]\) have overlapping interiors, but neither distinguished
    point (\(i\) or \(j\)) lies in the interior of the other interval.}

    We define two orderings of \(V(G_k)\) (the first will be used to avoid \(F_1\) and
    the second to avoid \(F_2\)).

    \paragraph{Ordering 1:} The vertices are ordered according to their second
    coordinate; i.e., \((a,i,b) < (c,j,d)\) if \(i < j\), with the vertices having the
    same second coordinate ordered arbitrarily.

    \paragraph{Ordering 2:} The vertices are ordered according to their first
    coordinate; i.e., \((a,i,b) < (c,j,d)\) if \(a < c\), with the vertices having the
    same first coordinate ordered arbitrarily.
\end{construction}

As we shall see, the graph \(G_k\) given by Construction \ref{const:F1,F2} is
\(F_1\)-free under Ordering 1 and \(F_2\)-free under Ordering 2. To obtain a sublinear
independence number, we will give weights to the vertices such that the weighted
independence number is sublinear. Blowing up the graph according to these vertex-weights
will provide the necessary construction.

When referring to the adjacency rule, we always mean the one in Construction
\ref{const:F1,F2}. We now show that \(G_k\) is \(F_1\)-free under Ordering 1 and
\(F_2\)-free under Ordering 2. Since we will later blow up the construction, we in fact
need to show that it is free of homomorphic images of \(F_i\) (where homomorphisms
\(\varphi\) are required to preserve the vertex-order, meaning that if \(x < y\) then
\(\varphi(x) \leq \varphi(y)\)).\footnote{The only homomorphic images of \(F_i\) are
\(F_i\) and \(K_3\).}
\begin{lem}\label{lem:F1-free, F2-free}
    For every \(k\ge3\), the ordered graph \(G_k\) has no homomorphic image of \(F_1\)
    under Ordering 1 and no homomorphic image of \(F_2\) under Ordering 2.
\end{lem}
\begin{proof}
    Consider first Ordering 1. Suppose, for contradiction, that there are vertices \(x_1
    \leq x_2 \leq x_3 \leq x_4\), \(x_r = (a_r,i_r,b_r)\), such that the map \(r \mapsto
    x_r\) is a homomorphism from \(F_1 = \{12,34,14\}\) to \(G_k\). By the definition of
    Ordering 1, we have \(i_1\le i_2\le i_3\le i_4\). Applying the adjacency rule (see
    Construction \ref{const:F1,F2}) to \(x_1x_2\) and \(x_3x_4\) gives \(b_1\le i_2\)
    and \(i_3\le a_4\), respectively. Since also \(i_2 \leq i_3\), we conclude that
    \(b_1 \leq a_4\). However, applying the adjacency rule to \(x_1x_4\) gives
    \(a_4<b_1\), a contradiction.

    Now consider Ordering 2. Suppose, for contradiction, that there are vertices \(x_1
    \leq x_2 \leq x_3 \leq x_4\), \(x_r = (a_r,i_r,b_r)\), such that the map \(r \mapsto
    x_r\) is a homomorphism from \(F_2\) to \(G_k\). By the definition of Ordering 2, we
    have \(a_1 \leq a_2 \leq a_3 \leq a_4\). Note that if two vertices \((a,i,b)\) and
    \((c,j,d)\) of \(G_k\) are adjacent and \(i<j\) then also \(a<c\), since \(a < i
    \leq c\). Hence, the fact that \(x_1,x_2\) are adjacent and that \(a_1 \leq a_2\)
    means that \(i_1 < i_2\) and hence \(i_1 \leq a_2 < b_1 \leq i_2\). Similarly, the
    fact that \(x_2,x_3\) are adjacent means that \(i_2 \leq a_3 < b_2 \leq i_3\). It
    follows that \(b_1 \leq i_2 \leq a_3 \leq a_4\). But then we cannot have the edge
    \(x_1x_4\), as this requires \(a_4 < b_1\) (again, because \(x_1x_4\) being adjacent
    and \(a_1 \leq a_4\) together imply that \(i_1 < i_4\), which in turn means that
    \(a_4 < b_1\) by the adjacency rule). This is a contradiction.
\end{proof}

Next, we will show that one can assign weights to the vertices of \(G_k\) so that the
weighted independence number is small.
\begin{thm}\label{thm:F1,F2 independence number}
    There exist integer weights \(w(a,i,b) \geq 1\) for \((a,i,b) \in \binom{[k]}{3}\)
    such that the sum of all weights is \(\Omega(k^3\log k)\), and such that the total
    weight of any independent set of \(G_k\) is at most \(2k^3\).
\end{thm}

Using Theorem \ref{thm:F1,F2 independence number}, we can easily deduce Theorem
\ref{thm:F1,F2}.
\begin{proof}[Proof of Theorem \ref{thm:F1,F2}]
    The upper bounds \(R_<(F_1,K_n) = O(n\log n)\) and \(R_<(F_2,K_n) = O(n\log n)\)
    follow from Lemma \ref{lem:e(I,J) = 1}. For the lower bounds, let \(G'_k\) be the
    blowup of \(G_k\) according to the weights \(w(a,i,b)\) given by Theorem
    \ref{thm:F1,F2 independence number}; namely, replace each vertex \((a,i,b)\) with an
    independent interval of \(w(a,i,b)\) vertices, preserve the order of these
    intervals, and replace every edge with a complete bipartite graph. Then for each
    \(i=1,2\), the graph \(G'_k\) has an \(F_i\)-free ordering by Lemma
    \ref{lem:F1-free, F2-free}. Furthermore, by Theorem
    \ref{thm:F1,F2 independence number}, we have \(|V(G'_k)| = \Omega(k^3\log k)\) and
    \(\alpha(G'_k) \leq 2k^3\). Indeed, the independence number of the blowup \(G'_k\)
    is equal to the maximum total weight of an independent set of \(G_k\). Finally, let
    \(n\) be sufficiently large and \(k:=\lfloor(\frac{n-1}{2})^{1/3}\rfloor\). Then
    \(\alpha(G'_k)\leq 2k^3 < n\) and \(k^3=\Theta(n)\). Hence, for \(i=1,2\),
    \[
        R_<(F_i,K_n)
        > |V(G'_k)|
        =
        \Omega(k^3\log k)
        =
        \Omega(n\log n).\qedhere
    \]
\end{proof}

\begin{proof}[Proof of Theorem \ref{thm:F1,F2 independence number}]
    For each \((a,i,b) \in \binom{[k]}{3}\) (where \(a < i < b\)), set
    \[
        w(a,i,b) :=
        \left\lfloor
        \frac{k^2}{(b-a)(b-a-1)}\right\rfloor\ge1.
    \]
    We first estimate the total sum of weights. For each \(2 \leq h \leq k-1\), there
    are exactly \((k-h)(h-1)\) triples \((a,i,b) \in \binom{[k]}{3}\) with \(b-a = h\).
    Furthermore, for each such triple \((a,i,b)\), we have \(w(a,i,b) \geq
    \frac{k^2}{h(h-1)}-1\). \nolinebreak Hence,
    \[
        \sum_{(a,i,b) \in \binom{[k]}{3}}w(a,i,b) \geq
        \sum_{h = 2}^{k-1} \left( \frac{k^2(k-h)}{h} - k^2 \right) \geq
        \frac{k^3}{2} \cdot \sum_{h = 2}^{\lfloor k/2 \rfloor} \frac{1}{h} - k^3 =
        \Omega(k^3\log k),
    \]
    using the partial sum of the harmonic series.

    Next, let \(I\subseteq V(G_k)\) be an independent set. For a pair \(1 \leq a < b
    \leq k\) with \(b-a \geq 2\), let
    \[
        m_{a,b}
        =\bigl|\{ i : a < i < b, \; (a,i,b) \in I\}\bigr|.
    \]
    Namely, \(m_{a,b}\) is the number of vertices of \(I\) of the form \((a,i,b)\).
    Define also
    \[
        \mathcal D := \{(a,b): m_{a,b} > 0\}.
    \]
    By the choice of the weights \(w(a,i,b)\), we have
    \begin{equation}
        \label{eq:F1F2-indep-bound}
        w(I) := \sum_{(a,i,b) \in I}w(a,i,b)
        \leq
        k^2\sum_{(a,b)\in\mathcal D}
        \frac{m_{a,b}}{(b-a)(b-a-1)} \; .
    \end{equation}
    It remains to show that the final sum is at most \(2k\), as this would imply that
    \(w(I) \leq 2k^3\), as required.

    We now arrive at the main idea of the proof. For \((a,b) \in \mathcal{D}\), define
    \begin{align*}
        L_a&=\max\{i-a: (a,i,b') \in I \text{ for some } b'\},\\
        R_b&=\max\{b-i: (a',i,b) \in I \text{ for some } a'\}.
    \end{align*}
    Note that both \(L_a\) and \(R_b\) are well-defined because, as \((a,b) \in
    \mathcal{D}\), there exists \(a < i < b\) such that \((a,i,b) \in I\). We now
    observe two crucial properties of the parameters \(m_{a,b}, L_a, R_b\). First, we
    have
    \begin{equation}
        \label{eq:m_ab <= L_a,R_b}
        m_{a,b}\leq\min\{L_a,R_b\},
    \end{equation}
    because every \(i\) satisfying \((a,i,b) \in I\) must belong to the interval
    \([a+1,a+L_a]\) and to the interval \([b-R_b,b-1]\). Second, we claim that
    \begin{equation}
        \label{eq:L_a or R_b is at most b-a}
        b-a>L_a
        \qquad \text{or} \qquad
        b-a>R_b.
    \end{equation}
    Indeed, suppose for a contradiction that \(L_a\geq b-a\) and \(R_b\geq b-a\). By the
    definition of \(L_a\) and \(R_b\), there exist \(j,i\) and \(a',b'\) such that
    \((a,j,b'),(a',i,b) \in I\) and \(j = a + L_a \geq b\), \(i = b - R_b \leq a\). But
    then \(i \leq a < b \leq j\), meaning that \((a',i,b),(a,j,b')\) are adjacent by the
    adjacency rule (in Construction \ref{const:F1,F2}). This contradicts that \(I\) is
    independent, proving \eqref{eq:L_a or R_b is at most b-a}.

    For each pair \((a,b) \in \mathcal{D}\), assign \((a,b)\) to \(b\) if \(b-a>R_b\)
    and to \(a\) otherwise; in the latter case, \eqref{eq:L_a or R_b is at most b-a}
    gives \(b-a>L_a\). Now, for each fixed \(b\), every pair \((a,b) \in \mathcal{D}\)
    assigned to \(b\) satisfies \(h := b-a>R_b\) and \(m_{a,b}\le R_b\) (by
    \eqref{eq:m_ab <= L_a,R_b}). Hence
    \[
        \sum_{\substack{(a,b)\in\mathcal D\\\text{assigned to }b}}
        \frac{m_{a,b}}{(b-a)(b-a-1)}
        \le R_b\sum_{h>R_b}\frac1{h(h-1)}
        =1.
    \]
    Here we replaced the finite sum over \(h > R_b\) with the infinite sum from
    \(R_b+1\) to infinity. Similarly, for each fixed \(a\), every pair \((a,b) \in
    \mathcal{D}\) assigned to \(a\) satisfies \(h := b-a>L_a\) and \(m_{a,b}\le L_a\)
    (by \eqref{eq:m_ab <= L_a,R_b}). Hence,
    \[
        \sum_{\substack{(a,b)\in\mathcal D\\\text{assigned to }a}}
        \frac{m_{a,b}}{(b-a)(b-a-1)}
        \le L_a\sum_{h>L_a}\frac1{h(h-1)}
        =1.
    \]
    Finally, summing over all \(k\) choices for \(b\) and all \(k\) choices for \(a\),
    we see that
    \[
        \sum_{(a,b) \in \mathcal{D}}\frac{m_{a,b}}{(b-a)(b-a-1)} \leq
        2k,
    \]
    and together with \eqref{eq:F1F2-indep-bound} we get \(w(I) \leq 2k^3\), as
    required.
\end{proof}

\subsection{Proof of Theorem \ref{thm:F3}}

Recall that \(F_3 = \{13, 24, 14\}\). We start by defining the ``twisted hypercube''
construction that is used to lower-bound \(R_<(F_3, K_n)\).

\begin{construction}\label{const:F3}
    Fix an integer \(h\geq 1\). We construct a random ordered graph \(G=G(h)\) on the
    vertex set \(\{0,1\}^h\), ordered lexicographically, as follows. For \(0\leq d\leq
    h\) and a string \(\tau=(\tau_1,\ldots,\tau_d)\in\{0,1\}^d\), let
    \[
        V_\tau :=
        \{x\in\{0,1\}^h:(x_1,\ldots,x_d)=\tau\}.
    \]
    Thus \(V_\varnothing=\{0,1\}^h\). For each \(\sigma\in\{0,1\}^{<h}\), choose
    independently and uniformly a perfect matching \(M_\sigma\) between \(V_{\sigma0}\)
    and \(V_{\sigma1}\).\footnote{Here, for \(i=0,1\), we use \(\sigma i\) to denote the
    string obtained from \(\sigma\) by appending \(i\).} The edge set of \(G\) is the
    union of these matchings.
\end{construction}

Construction \ref{const:F3} can be described inductively: place a uniformly random
perfect matching between the first \(2^{h-1}\) vertices and the last \(2^{h-1}\)
vertices, and then iterate independently inside each half.

Recall that \(F_3 = \{13,24,14\}\). We now show that the above construction is always
\(F_3\)-free.

\begin{lem}\label{lem:F3-free}
    Every realization of the graph \(G\) given by Construction \ref{const:F3} is
    \(F_3\)-free.
\end{lem}

\begin{proof}
    Suppose by contradiction that \(x_1<x_2<x_3<x_4\) satisfy
    \[
        x_1x_3, \; x_1x_4, \; x_2x_4\in E(G).
    \]
    Let \(\sigma \in \{0, 1\}^{< h}\) be the longest common prefix of \(x_1\) and
    \(x_4\) (considered as strings in \(\{0,1\}^h\)). Then
    \[
        x_1\in V_{\sigma0},
        \qquad
        x_4\in V_{\sigma1},
        \qquad
        x_1x_4\in M_\sigma.
    \]
    In particular, \(x_1,x_4 \in V_{\sigma}\) (and \(x_1x_4 \in M_{\sigma}\)). As
    \(V_\sigma\) is an interval in the vertex order, we have \(x_2,x_3\in V_\sigma\). If
    \(x_2\in V_{\sigma0}\), then \(x_2x_4\in M_\sigma\), contradicting that \(M_\sigma\)
    is a matching. Otherwise \(x_2\in V_{\sigma1}\), and hence \(x_3\in V_{\sigma1}\).
    It follows that \(x_1x_3\in M_\sigma\), again contradicting that \(M_\sigma\) is a
    matching.
\end{proof}

Next, we bound the independence number of the construction.
\begin{thm}\label{thm:F3 independence number}
    Let \(G = G(h)\) be given by Construction \ref{const:F3}. Then \(\alpha(G) =
    O(\frac{2^h \log h}{h})\) w.h.p.
\end{thm}

Using Theorem \ref{thm:F3 independence number}, we can now easily deduce Theorem
\ref{thm:F3}.

\begin{proof}[Proof of Theorem \ref{thm:F3}]
    For the lower bound, let \(n\) be sufficiently large and set
    \(h:=\lfloor\frac{1}{2}\log_2 n\rfloor\). By Lemma \ref{lem:F3-free} and Theorem
    \ref{thm:F3 independence number}, choose an \(F_3\)-free outcome \(H\) of \(G(h)\)
    satisfying \(\alpha(H)\leq C2^h\log h/h\) for some constant \(C>0\). Let
    \(\widetilde{H}\) be the disjoint union of \(\lfloor(n-1)/\alpha(H)\rfloor\) copies
    of \(H\). Then \(\widetilde{H}\) is \(F_3\)-free and \(\alpha(\widetilde{H})<n\).
    Since \(\alpha(H)\leq 2^h\leq\sqrt n=o(n)\), we have
    \(\lfloor(n-1)/\alpha(H)\rfloor\geq n/(2\alpha(H))\). Therefore,
    \[
        R_<(F_3,K_n)
        > |V(\widetilde{H})|
        =
        \left\lfloor
        \frac{n-1}{\alpha(H)}
        \right\rfloor 2^h
        \geq
        \frac{n2^h}{2\alpha(H)}
        \geq
        \frac{nh}{2C\log h}
        =
        \Omega\left(
        \frac{n\log n}{\log\log n}
        \right).
    \]

    For the upper bound, we use the fact that \(\ex_<(n,F_3) = \Theta(n\log n)\), see
    \cite{PachTardos}. Thus, every \(n\)-vertex \(F_3\)-free graph has average degree
    \(d = O(\log n)\). We now additionally claim that in an \(F_3\)-free graph, the
    subgraph induced by the neighborhood of any vertex \(v\) is 4-colorable. Consider
    first the forward neighborhood of \(v\). If there are forward neighbors \(x<y<z\) of
    \(v\) such that \(xz \in E(G)\), then the edges \(vy,xz,vz\) form a copy of \(F_3\),
    a contradiction. It follows that there can only be edges between consecutive forward
    neighbors of \(v\). Thus, the subgraph induced by the forward neighborhood of \(v\)
    is contained in a path, and hence has chromatic number at most 2. The same argument
    applies to the backward neighborhood of \(v\). Thus, \(\chi(G[N(v)]) \leq 4\).

    A result of Alon \cite{Alon} states that for every \(r \geq 2\), if \(G\) is an
    \(n\)-vertex graph with average degree \(d\) and every neighborhood in \(G\) is
    \(r\)-colorable, then \(\alpha(G) = \Omega_r(\frac{n \log d}{d})\). Applying this
    with \(r=4\), we conclude that every \(n\)-vertex \(F_3\)-free ordered graph \(G\)
    has \(\alpha(G) = \Omega(\frac{n \log\log n}{\log n} )\). This proves that
    \(R_<(F_3,K_n) =O (\frac{n\log n}{\log\log n}) \).
\end{proof}

It remains to prove Theorem \ref{thm:F3 independence number}. We will use the following
lemma, which bounds the probability that a random injection misses a prescribed subset.

\begin{lem}\label{lem:matching}
    Let \(A \neq \varnothing\) and \(B\) be two sets, let \(S \subseteq A\), and let
    \(\phi : B \rightarrow A\) be a uniformly-random injection. Then
    \[
        \Pr\bigl(\phi(B) \cap S=\varnothing\bigr)
        \le \exp\!\left(-\frac{|B||S|}{|A|}\right).
    \]
\end{lem}
\begin{proof}
    Set \(a=|A|\), \(b=|B|\) and \(s = |S|\). If \(b+s>a\), the probability is zero.
    Otherwise,
    \[
        \Pr\bigl(\phi(B) \cap S=\varnothing\bigr)
        =\frac{\binom{a-s}{b}}{\binom{a}{b}}
        \le \left(1-\frac{s}{a}\right)^b
        \le e^{-bs/a}.\qedhere
    \]
\end{proof}

\begin{cor}\label{cor:matching}
    Let \(M\) be a uniformly-random perfect matching between two sets \(A\) and \(B\) of
    size \(r \geq 1\). Let \(S \subseteq A\) and \(T \subseteq B\). Then
    \[
        \Pr\bigl(M\cap(S \times T)=\varnothing\bigr)
        \le \exp\!\left(-\frac{|S||T|}{r}\right).
    \]
\end{cor}
\begin{proof}
    The restriction \(M|_T\) is a random injection from \(T\) to \(A\).
\end{proof}

\begin{proof}[Proof of Theorem \ref{thm:F3 independence number}]
    Put \(n := |V(G)| = 2^h\). Fix \(S\subseteq V(G)\) of size \(s\), and for each
    string \(\tau \in \{0,1\}^{\leq h}\), set
    \[
        s_\tau := |S\cap V_\tau|.
    \]
    Recall that \(E(G)=\bigcup_{\sigma\in\{0,1\}^{<h}} M_\sigma\), and that
    \(M_{\sigma}\) is a random matching between the two sets \(V_{\sigma 0},V_{\sigma
    1}\) of size \(2^{h-|\sigma|-1}\) each. Also, these matchings are chosen
    independently. Hence, Corollary~\ref{cor:matching} yields
    \begin{equation}\label{eq:F3 independence 1}
        \mathbb P(S\text{ is independent})
        = \prod_{\sigma \in \{0, 1\}^{<h}} \mathbb P\left(M_\sigma \cap \binom{S}{2} =
            \varnothing\right)
        \le \exp\left(-\sum_\sigma
        \frac{s_{\sigma0}s_{\sigma1}}{2^{h-|\sigma|-1}}\right).
    \end{equation}
    Note that \(s_{\sigma} = s_{\sigma0} + s_{\sigma1}\) and therefore
    \(2s_{\sigma0}s_{\sigma1} = s_\sigma^2 - s_{\sigma0}^2 - s_{\sigma1}^2\). Hence, we
    have
    \begin{align}
        \label{eq:F3 independence 2}
        \sum_{\sigma\in\{0,1\}^{<h}} \frac{s_{\sigma0}s_{\sigma1}}{2^{h-|\sigma|-1}}
        &= \nonumber \sum_{\ell=0}^{h-1}\frac{2^\ell}{2^h} \left(
        \sum_{\sigma\in\{0,1\}^\ell}s_\sigma^2
        - \sum_{\tau\in\{0,1\}^{\ell+1}}s_\tau^2 \right) \\
        &= \nonumber \frac{s_\varnothing^2}{2^h} + \sum_{\ell=1}^{h-1} \left(
            \frac{2^\ell -
        2^{\ell-1}}{2^h} \right) \left( \sum_{\tau\in\{0,1\}^\ell} s_\tau^2 \right) -
        \frac{2^{h-1}}{2^h} \sum_{\tau\in\{0,1\}^h} s_\tau^2 \\
        &= \frac{s^2}{2^h} + \frac{1}{2^{h+1}}\sum_{\ell=1}^{h-1}2^\ell
        \left(\sum_{\tau\in\{0,1\}^\ell}s_\tau^2\right) - \frac{s}{2}.
    \end{align}
    Here, in the last equality we used that \(s_\varnothing = s\) (as \(V_{\varnothing}
    = \{0,1\}^h\)), and that \(\sum_{\tau\in\{0,1\}^h}s_\tau^2 = s\), as \(s_\tau \in
    \{0, 1\}\) when \(\tau \in \{0, 1\}^h\). Finally, Cauchy--Schwarz gives
    \[
        2^\ell \left(\sum_{\tau\in\{0,1\}^\ell}s_\tau^2\right) =
        \left(\sum_{\tau\in\{0,1\}^\ell}1^2\right)\left(\sum_{\tau\in\{0,1\}^\ell}s_\tau^2\right)\ge
        \left(\sum_{\tau\in\{0,1\}^\ell}s_\tau\right)^2 = s^2
    \]
    for each \(1\le\ell\leq h-1\), using that the sets \((V_{\tau} : \tau \in
    \{0,1\}^{\ell})\) partition the vertex set. Plugging this into
    \eqref{eq:F3 independence 2} and then into \eqref{eq:F3 independence 1}, we obtain
    \[
        \mathbb P(S\text{ is independent}) \le \exp\left( -\frac{s^2}{2^h} -
        \frac{(h-1)s^2}{2^{h + 1}}
        + \frac{s}{2} \right) = \exp\left( -\frac{(h+1)s^2}{2^{h + 1}} + \frac{s}{2}
            \right).
    \]
    Set \(s := \lceil \frac{2^{h+1}\log h}{h} \rceil\). The union bound over all
    \(\binom{n}{s} = \binom{2^h}{s}\) choices for \(S\) yields
    \begin{align*}
        \mathbb P\bigl(\alpha(G) \ge s\bigr)
        &\le \binom{2^h}{s}
        \exp\left(-\frac{(h+1)s^2}{2^{h+1}} + \frac{s}{2}\right) \\
        &\le \left(
        \frac{e 2^h}{s}
        \exp\left(-\frac{(h+1)s}{2^{h+1}} + \frac{1}{2}\right)
        \right)^s = o_h(1).
    \end{align*}
    Here, the last equality holds since \(\frac{2^h}{s} = O(\frac{h}{\log h})\) and
    \[
        \exp\left(- \frac{(h+1)s}{2^{h+1}} + \frac{1}{2}\right) \leq
        e^{-\log h + \frac{1}{2}} = O\left(\frac{1}{h}\right).
    \]
    The theorem follows.
\end{proof}

\subsection{Proof of Theorem \ref{thm:F4}}

Recall that \(F_4 = \{12, 13, 24\}\). We use the following construction to lower-bound
\(R_<(F_4, K_n)\).
\begin{construction}\label{const: Gh}
    For each \(h \geq 0\), we construct an ordered graph \(G_h\) and an independent set
    \(A_h \subseteq V(G_h)\), such that \(A_h\) is a prefix of the vertex-order of
    \(G_h\). The graph \(G_0\) is the 1-vertex graph, and \(A_0 := V(G_0)\).

    Suppose that \((G_h,A_h)\) has been constructed. To construct \((G_{h+1},A_{h+1})\),
    take three disjoint copies \(H_1,H_2,H_3\) of \(G_h\), and for \(x \in A_h\), let
    \(x_i\) denote the copy of \(x\) in \(H_i\). Let
    \[
        A_{h+1} := \{(i,x) : i \in [3],\; x \in A_h\}
    \]
    be a new set, disjoint from \(H_1,H_2,H_3\), consisting of three disjoint copies of
    \(A_h\). The vertex set of \(G_{h+1}\) is \(A_{h+1} \cup V(H_1) \cup V(H_2) \cup
    V(H_3)\), ordered so that
    \[
        A_{h+1} < V(H_1) < V(H_2) < V(H_3),
    \]
    with the vertices of \(A_{h+1}\) ordered arbitrarily and each \(H_i\) retaining its
    original order. Clearly, \(A_{h+1}\) is a prefix of \(V(G_{h+1})\). The edges of
    \(G_{h+1}\) consist of the edges of \(H_1,H_2,H_3\) and the following edges: For
    each \(i \in [3]\) and \(x \in A_h\), join \((i,x) \in A_{h+1}\) to \(x_i\) and to
    every vertex of \(N_{H_j}(x_j)\) for \(j \in [3] \setminus \{i\}\).
\end{construction}

Let us record the sizes of \(V(G_h)\) and \(A_h\).
\begin{claim}\label{claim:F4 construction size}
    For every \(h \geq 0\), it holds that \(|V(G_h)| = (h+1)3^h\) and \(|A_h| = 3^h\).
\end{claim}
\begin{proofclaim}
    By induction on \(h\). The case \(h=0\) is trivial. For the induction step,
    \(|A_{h+1}| = 3|A_h| = 3^{h+1}\), and
    \[
        |V(G_{h+1})| = 3(|V(G_h)| + |A_h|) = 3 \cdot (h+2)3^h = (h+2)3^{h+1},\]
        as required.
    \end{proofclaim}

    Next, we prove that this construction is \(F_4\)-free.
    \begin{lem}\label{lem:F4-free}
        The graph \(G_h\) is \(F_4\)-free.
    \end{lem}
    \begin{proof}
        We argue by induction on \(h\). The case \(h=0\) is trivial. Suppose that the
            lemma holds
        for \(h\), and suppose for a contradiction that \(a<b<c<d\) form a copy of
            \(F_4\) in
        \(G_{h+1}\), so that \(ab,ac,bd \in E(G_{h+1})\).

        Recall that there are no edges between the blocks \(H_1,H_2,H_3\), and that
            \(A_{h+1}\)
        precedes them in the vertex order. Hence, if \(a\) lies in some block \(H_i\),
            then so do
        \(b,c,d\), since \(F_4\) is connected. But \(H_i \cong G_h\) is \(F_4\)-free by
            the
        induction hypothesis.

        Thus, \(a \in A_{h+1}\), and we may write
        \(a=(i,x)\) for some \(i \in [3]\) and \(x \in A_h\).
        Since \(A_{h+1}\) is independent, \(b \notin A_{h+1}\). The
        vertices \(b\) and \(d\) then lie in the same block \(H_j\), since \(bd\) is an
            edge. As
        \(H_j\) is an interval and \(b<c<d\), the vertex \(c\) also lies in \(H_j\).
            Thus, \((i,x)\)
        has two neighbors in \(H_j\), namely \(b\) and \(c\). By the definition of
            \(G_{h+1}\), this is only possible if \(j \neq i\)
        and \(b,c \in N_{H_j}(x_j)\) (indeed, if \(j=i\) then \((i,x)\) has only one
            neighbor in \(H_j\)). Since \(x_j\) belongs to the copy of \(A_h\) in
            \(H_j\), it
        precedes all its neighbors. Hence, \(x_j<b<c<d\) form a copy of \(F_4\) in
            \(H_j\), a
        contradiction.
    \end{proof}

    To handle the independence number, we prove the following:
    \begin{thm}\label{thm:F_4 independence number}
        For every \(h \geq 0\), the graph \(G_h\) given by
        Construction \ref{const: Gh} satisfies \( \alpha(G_h) \leq 2 \cdot 3^h. \)
    \end{thm}

    \begin{proof}
        The theorem follows from the stronger formulation below, which we prove by
            induction.
        \begin{claim}\label{claim:F4 independence number}
            If \(I\) is independent in \(G_h\),
            then\footnote{Here,
            \(N(I) := \{v \in V(G_h) : uv \in E(G_h) \text{ for some } u \in I\}\) is
                the
            neighborhood of the set \(I\).}
            \[
            |I \setminus A_h| \leq 2|N(I) \cap A_h|.
        \]
    \end{claim}

    Claim \ref{claim:F4 independence number} implies that
    \[
        |I| =
        |I \cap A_h| + |I \setminus A_h|
        \leq |I \cap A_h| + 2|N(I) \cap A_h|
        \leq 2|A_h|
        = 2 \cdot 3^h,
    \]
    proving Theorem \ref{thm:F_4 independence number}. Here, in the last inequality we
    used that \(I \cap N(I) = \varnothing\) (as \(I\) is independent).

    \begin{proofclaim}[Proof of Claim \ref{claim:F4 independence number}]
        We argue by induction on \(h\). The case \(h=0\) is trivial. Suppose that the
        claim holds for \(h\), and let \(I \subseteq V(G_{h+1})\) be independent. For
        \(i \in [3]\), let \(I_i := I \cap V(H_i)\), and let \(B_i\) be the copy of
        \(A_h\) in \(H_i\). Then, the induction hypothesis gives
        \begin{equation}
            \label{eq:sum over copies}
            |I \setminus A_{h+1}|
            = \sum_{i=1}^{3} |I_i|
            =
            \sum_{i=1}^{3}
            \bigl(|I_i \cap B_i| + |I_i \setminus B_i| \bigr)
            \leq \sum_{i=1}^{3}
            \bigl(|I_i \cap B_i| + 2|N_{H_i}(I_i) \cap B_i|\bigr).
        \end{equation}
        Since each \(B_i\) is a copy of \(A_h\), we can rewrite the RHS of
        \eqref{eq:sum over copies} as a sum over \(x \in A_h\), as follows. For \(x \in
        A_h\), let \(P_x\) be the set of \(i \in [3]\) with \(x_i \in I_i\), let \(Q_x\)
        be the set of \(i \in [3]\) with \(x_i \in N_{H_i}(I_i)\), and let \(R_x\) be
        the set of \(i \in [3]\) with \((i,x) \in N(I)\). Then the RHS of
        \eqref{eq:sum over copies} is \(\sum_{x \in A_h} (|P_x|+2|Q_x|)\), while \(
        |N(I) \cap A_{h+1}| = \sum_{x \in A_h} |R_x|. \) Thus, to complete the induction
        step for \(h+1\), it suffices to show that
        \begin{equation}
            \label{eq:local}
            |P_x| + 2|Q_x| \leq 2|R_x|
        \end{equation}
        for every \(x \in A_h\). Fix \(x \in A_h\). We use three facts:
        \begin{enumerate}[label=(\roman*)]
            \item \(P_x\) and \(Q_x\) are disjoint because \(I\) is independent, so
            \(|P_x|+|Q_x| \leq 3\);

            \item \(P_x \subseteq R_x\);
            \item For every \(i \in Q_x\) it holds that \([3]\setminus\{i\} \subseteq
            R_x\).

        \end{enumerate}
        Item (i) is immediate. Item (ii) holds by the definition of \(G_{h+1}\): For
        every \(x \in A_h\), the vertices \((i,x)\) and \(x_i\) are adjacent. Hence, if
        \(x_i \in I_i\) then \((i,x) \in N(I)\), implying that \(i \in R_x\). We now
        prove Item (iii). So suppose that \(i \in Q_x\), and let \(y_i \in I_i\) such
        that \(x_i,y_i\) are adjacent. By the definition of \(G_{h+1}\), we then have
        that \((j,x),y_i\) are adjacent for all \(j \in [3] \setminus \{i\}\).
        Therefore, for every \(j \in [3] \setminus \{i\}\) it holds that \((j,x) \in
        N(I)\) and hence \(j \in R_x\), as required.

        We now prove \eqref{eq:local} in three cases.
        \begin{enumerate}
            \item \(|Q_x| \geq 2\). By (iii), \(R_x=[3]\), so (i) gives
            \[
                |P_x|+2|Q_x| \leq 3+|Q_x| \leq 6=2|R_x|.
            \]

            \item \(|Q_x|=1\). This follows since \(|R_x| \geq 2\) by (iii) and \(|P_x|
            \leq 2\) by (i) so
            \[
                |P_x|+2|Q_x| \leq 2+2 \leq 2|R_x|.
            \]

            \item \(Q_x=\emptyset\). Then \eqref{eq:local} follows directly from
            (ii).\qedhere
        \end{enumerate}
    \end{proofclaim}
\end{proof}

\begin{proof}[Proof of Theorem \ref{thm:F4}]
    The upper bound follows from Lemma \ref{lem:e(I,J) = 1}. For the lower bound, let
    \(n\) be sufficiently large and set \(h:=\lfloor\frac{1}{2}\log_3 n\rfloor\). By
    Lemma \ref{lem:F4-free} and Theorem \ref{thm:F_4 independence number}, the graph
    \(G_h\) is \(F_4\)-free and satisfies \(|V(G_h)|=(h+1)3^h\) and \(\alpha(G_h)\leq
    2\cdot 3^h\). Let \(\widetilde G\) be the disjoint union of
    \(\lfloor(n-1)/\alpha(G_h)\rfloor\) copies of \(G_h\). Then \(\widetilde G\) is
    \(F_4\)-free (since \(F_4\) is connected), and \(\alpha(\widetilde G)<n\). As
    \(\alpha(G_h)\leq 2\sqrt n=o(n)\), we have \(\lfloor(n-1)/\alpha(G_h)\rfloor\geq
    n/(2\alpha(G_h))\). Therefore,
    \[
        R_<(F_4,K_n)
        > |V(\widetilde G)|
        \geq \frac{n|V(G_h)|}{2\alpha(G_h)}
        \geq \frac{n(h+1)}{4}
        = \Omega(n\log n).\qedhere
    \]
\end{proof}

\subsection{Proof of Theorem \ref{thm:R}}

Recall that \(R = \{13, 23, 14, 25\}\). We use the following construction to lower-bound
\(R_<(R, K_n)\).

\begin{construction}\label{const:randomized-Tutte}
    Fix an integer \(q\geq1\). We construct random ordered graphs \(G=G(q,h)\), \(1\leq
    h\leq q\),\footnote{The assumption \(h \leq q\) ensures that the sets \(B_i\) are
    not larger than \(A\), as we need an injection from \(B_i\) to \(A\).} where
    \(G(q,h)\) has \(h q^{h-1}\) vertices. The graph \(G(q,1)\) is simply the 1-vertex
    graph. For \(h\geq2\), take consecutive intervals \(A<B_1<\dots<B_q\), where
    \[
        |A|=q^{h-1},\qquad |B_i|=(h-1)q^{h-2}.
    \]
    Place mutually independent copies of \(G(q,h-1)\) on \(B_1,\dots,B_q\), declare
    \(A\) to be independent, and place no edges between \(B_i\) and \(B_j\) for \(i\neq
    j\). Independently for every \(i\in[q]\), choose a uniformly random injection
    \(\phi_i\colon B_i\to A\) and add the edges \(b\phi_i(b)\) for all \(b\in B_i\).
\end{construction}

Note that the number of vertices of \(G(q,h)\) is indeed \( |A| + |B_1| + \dots + |B_q|
= \linebreak q^{h-1} + q \cdot (h-1) q^{h-2} = hq^{h-1}\), as required.

\begin{lem}\label{lem:K-free}
    Every realization of \(G(q,h)\) is \(R\)-free.
\end{lem}

\begin{proof}
    We argue by induction on \(h\). The case \(h=1\) is trivial. Suppose that \(h\geq2\)
    and that \(v_1<v_2<v_3<v_4<v_5\) form a copy of \(R\), so that
    \(v_1v_3,v_2v_3,v_1v_4,v_2v_5\in E(G)\). If \(v_1\in B_i\) (for some \(i \in [q]\)),
    then the entire copy lies in \(B_i\), since \(A\) precedes \(B_i\) and there are no
    edges between the sets \(B_1,\dots,B_q\). This contradicts the induction hypothesis.

    Thus, \(v_1\in A\). Since \(A\) is independent, \(v_3\notin A\). Also, \(v_2\notin
    A\), because \(v_3\) has (exactly) one neighbor in \(A\), by construction. The edges
    \(v_2v_3,v_2v_5\) then imply that \(v_2,v_3,v_5\) lie in the same block \(B_i\).
    Since \(B_i\) is an interval, \(v_4\) also lies in \(B_i\). This is impossible
    because \(v_1\) has at most one neighbor in \(B_i\).
\end{proof}

To prove the lower bound in Theorem \ref{thm:R}, we will take \(h = q\) in Construction
\ref{const:randomized-Tutte}. The following theorem handles the independence number.

\begin{thm}\label{thm:randomized-Tutte-independence}
    Let \(G=G(q,q)\) be given by Construction \ref{const:randomized-Tutte} and
    \(n:=|V(G)| = q^q\). Then w.h.p.,
    \[
        \alpha(G)
        =O\left(\frac{n\log q}{q}\right)
        =O\left(\frac{n(\log\log n)^2}{\log n}\right).
    \]
\end{thm}

\begin{proof}
    The main step in the proof is to show that for every \(1\leq h\leq q\) and every
    \(S\subseteq V(G(q,h))\) of size \(s\),
    \begin{equation}
        \label{eq:Tutte independent set}
        \mathbb P(S\text{ is independent})
        \leq\exp\left(\frac{s}{2}-\frac{s^2}{2q^{h-1}}\right).
    \end{equation}
    The proof is by induction on \(h\). The base case \(h=1\) is trivial because
    \(G(q,1)\) is the 1-vertex graph. Now let \(h\geq2\), and set
    \[
        m:=|A|=q^{h-1},\qquad
        a:=|S\cap A|,\qquad
        b_i:=|S\cap B_i|,\qquad
        b:=\sum_{i=1}^q b_i=s-a.
    \]
    By Lemma \ref{lem:matching}, the image \(\phi_i(S\cap B_i)\) is disjoint from
    \(S\cap A\) with probability at most \(e^{-ab_i/m}\). Since the injections
    \(\phi_1,\dots,\phi_q\) are independent, the probability that for every \(i \in
    [q]\) there are no edges inside \(S\) between \(A\) and \(B_i\) is at most
    \(e^{-a(b_1+\dots+b_q)/m} = e^{-ab/m}\). Additionally, by the induction hypothesis,
    the probability that \(S \cap B_i\) is independent for every \(i \in [q]\) is at
    most
    \[
        \prod_{i=1}^q \exp\left(\frac{b_i}{2}-\frac{b_i^2}{2q^{h-2}}\right) =
        \exp\left( \frac b2 -\frac{q}{2m}\sum_{i=1}^q b_i^2 \right).
    \]
    Therefore,
    \[
        \begin{aligned}
            \mathbb P(S\text{ is independent})
            &\leq \exp\left(
            -\frac{ab}{m} + \frac{b}{2}
            - \frac{q}{2m}\sum_{i=1}^q b_i^2
            \right) \\
            &\leq \exp\left(\frac{b}{2} - \frac{ab}{m} - \frac{b^2}{2m}\right)
            \leq \exp\left(\frac{s}{2} - \frac{s^2}{2m}\right).
        \end{aligned}
    \]
    Here, the second inequality uses Cauchy--Schwarz. As for the last inequality, upon
    plugging in \(s=a+b\) and simplifying, this inequality becomes \(a(m-a) \geq 0\),
    which holds as \(a \leq m\). This proves \eqref{eq:Tutte independent set}.

    Now set \(h:=q\) and \(s:=\lceil4n\log q/q\rceil\). By
    \eqref{eq:Tutte independent set} and the union bound, and using that \(n = q^q\), we
    get
    \[
        \mathbb P\bigl(\alpha(G)\geq s\bigr)\leq\binom ns\exp\left(\frac
        s2-\frac{qs^2}{2n}\right)\leq\left(\frac{en}{s}
        \exp\left(\frac12-\frac{qs}{2n}\right)\right)^s
        \leq\left(\frac{e^{3/2}}{4q\log q}\right)^s=o(1).
    \]
    Finally, as \(s=O(n\log q/q)\) and \(q=\Theta(\log n/\log\log n)\), the theorem
    follows.
\end{proof}

\begin{proof}[Proof of Theorem \ref{thm:R}]
    For the lower bound, let \(n\) be sufficiently large and set
    \[
        q:=\lfloor\log n/(2\log\log n)\rfloor.
    \]
    By Lemma \ref{lem:K-free} and Theorem \ref{thm:randomized-Tutte-independence},
    choose an \(R\)-free outcome \(H\) of \(G(q,q)\) satisfying \(\alpha(H)=O(q^q\log
    q/q)\). Let \(\widetilde H\) be the disjoint union of
    \(\lfloor(n-1)/\alpha(H)\rfloor\) copies of \(H\). Since \(R\) is connected,
    \(\widetilde H\) is \(R\)-free, and \(\alpha(\widetilde H)<n\). Since
    \(\alpha(H)\leq q^q\leq\sqrt n=o(n)\), we have \(\lfloor(n-1)/\alpha(H)\rfloor\geq
    n/(2\alpha(H))\). Therefore,
    \[
        R_<(R,K_n)>|V(\widetilde H)|
        =\left\lfloor\frac{n-1}{\alpha(H)}\right\rfloor q^q
        \geq\frac{nq^q}{2\alpha(H)}=\Omega\left(\frac{nq}{\log q}\right)
        =\Omega\left(\frac{n\log n}{(\log\log n)^2}\right).
    \]

    For the upper bound, we proceed similarly as in the proof of Theorem \ref{thm:F3}.
    It is known that \(\ex_<(n,R) = \Theta(n\log n)\); the upper bound follows e.g.~from
    \cite[Lemma 2]{PachTardos} (and for the lower bound see \cite{FKMV}). Hence, every
    \(n\)-vertex \(R\)-free ordered graph \(G\) has average degree \(O(\log n)\). We
    claim that additionally, for every \(v \in V(G)\), the neighborhood \(G[N(v)]\) has
    bounded chromatic number (namely, bounded by an absolute constant). The reason is
    that both subgraphs \(R - \{1\}\) and \(R - \{5\}\) (obtained from \(R\) by removing
    the first or last vertex) are \(\chi\)-bounding (this is easily seen to follow from
    Theorem \ref{thm:charac}). The forward neighborhood of \(v\) avoids \(R - \{1\}\),
    and the backward neighborhood of \(v\) avoids \(R - \{5\}\). This proves our claim
    that \(\chi(G[N(v)])\) is bounded by some absolute constant \(C\). By the result of
    Alon \cite{Alon}, it follows that \(\alpha(G) \geq \Omega_C(\frac{n\log\log n}{\log
    n})\), which in turn shows that \(R_<(R,K_n) = O (\frac{n\log n}{\log\log n})\), as
    required.
\end{proof}

\section{Lower bounds via Davenport-Schinzel sequences}\label{sec:DS}
A \defn{Davenport-Schinzel} (\defn{DS} for short) \defn{sequence of order 3} is a
sequence \(\lambda_1,\dots,\lambda_N\) over a given alphabet such that no two
consecutive symbols are equal and such that the sequence is \(ababa\)-free, meaning that
there are no two distinct symbols \(a,b\) and indices \(1 \leq i_1 < \dots < i_5 \leq
N\) with \(\lambda_{i_1} = \lambda_{i_3} = \lambda_{i_5} = a\) and \(\lambda_{i_2} =
\lambda_{i_4} = b\). It is well-known (see \cite{HartSharir,Nivasch}) that there exist
DS-sequences of order 3 over the alphabet \([n]\) having length \(N = \Theta(n
\alpha(n))\) (where \(\alpha(n)\) is the inverse Ackermann function); in fact, Nivasch's
construction gives such a sequence of length \(N = (2-o(1))n \cdot \alpha(n)\) (and the
constant \(2\) is best possible). We need the following fact:
\begin{lem}\label{lem:splitting DS sequence}
    Every order-3 DS-sequence \(\lambda_1,\dots,\lambda_N\) over the alphabet \([n]\)
    can be split into at most \(2n\) intervals such that each interval consists of
    distinct symbols.
\end{lem}
\begin{proof}
    For each symbol \(a \in [n]\), mark the first and last appearances of \(a\) in the
    sequence, and split the sequence at these indices. More precisely, if \(1 = i_1 <
    \dots < i_m = N\) are the marked indices (note that \(i_1 = 1\) and \(i_m = N\)
    because the index \(1\) is the first appearance of the symbol \(\lambda_1\) and the
    index \(N\) is the last appearance of the symbol \(\lambda_N\)), then split \([N]\)
    into the intervals \([1=i_1,i_2),[i_2,i_3),\dots,[i_{m-1},i_m = N), \{N\}\). This
    splits the sequence into \(m \leq 2n\) intervals. Suppose by contradiction that some
    symbol \(b\) appears twice in some interval \(I\). Since consecutive symbols are
    distinct, there is some other symbol \(a\) between these two appearances of \(b\).
    This appearance of \(a\) is not one of the endpoints of the interval \(I\), and so
    this is neither the first nor the last appearance of \(a\). This means that there is
    an appearance of \(a\) to the left of \(I\) and an appearance of \(a\) to the right
    of \(I\) (possibly including the endpoints of \(I\)). This gives an
    \(ababa\)-pattern, a contradiction.
\end{proof}

Take an order-3 DS-sequence of length \(N = \Theta(n \alpha(n))\), and partition it into
at most \(2n\) intervals consisting of distinct symbols, as in Lemma
\ref{lem:splitting DS sequence}. Setting \(h := \lfloor \frac{N}{4n} \rfloor\),
partition each such interval into intervals of length \(h\) and a leftover of length
less than \(h\). The total size of all leftovers is less than \(h \cdot 2n \leq
\frac{N}{2}\). Thus, we obtain at least \(\frac{N}{2h} \geq n\) intervals of length
\(h\), each consisting of distinct symbols.
\begin{fact}\label{fact:DS sequence}
    There exists an \(ababa\)-free sequence over the alphabet \([n]\) which can be
    partitioned into \(n\) intervals, each consisting of \(h = \Omega(\alpha(n))\)
    distinct symbols.
\end{fact}

We are now ready to describe our construction used to prove that \(R_<(F,K_n) =
\omega(n)\) for \(F = F_0 = \{13,15,24\}\) and \(F = \{13,35,24\}\). For a symbol
\(\lambda \in [n]\), we use \(f(\lambda)\) to denote the index of the first appearance
of \(\lambda\) in the sequence.

\begin{construction}\label{const:DS}
    Consider the \(ababa\)-free sequence given by Fact \ref{fact:DS sequence}, and let
    \(B_1 < \dots < B_n\) be its partition into intervals of size \(h\), each consisting
    of distinct symbols. For each \(i \in [n]\), sample an entry of \(B_i\) uniformly at
    random and independently, and let \(\lambda_i\) be the corresponding symbol.

    Define a graph \(G\) on \([n]\) as follows: for \(i < j\), let \(ij \in E(G)\) if
    and only if the following hold:
    \begin{enumerate}
        \item \(\lambda_i \in B_j\) (namely, \(\lambda_i\) appears in the interval
        \(B_j\)).
        \item For every \(i < k < j\), \(\lambda_k \neq \lambda_i\).
        \item \(f(\lambda_i) > f(\lambda_j)\).
    \end{enumerate}
\end{construction}

The ordered graph given by Construction \ref{const:DS} is \(\{13,15,24\}\)- and
\(\{13,35,24\}\)-free with probability 1.
\begin{lem}\label{lem:DS construction freeness}
    Every realization of the ordered graph \(G\) given by Construction \ref{const:DS} is
    \(\{13,15,24\}\)-free and \(\{13,35,24\}\)-free.\footnote{We note that for
    \(\{13,15,24\}\)-freeness we may drop Condition 3 in Construction \ref{const:DS},
    and for \(\{13,35,24\}\)-freeness we may drop Condition 2. We opted for presenting
    one construction which works for both ordered forests.}
\end{lem}
\begin{proof}
    We first prove that the construction is \(\{13,15,24\}\)-free. Suppose by
    contradiction that there is a copy of \(\{13,15,24\}\) on the vertices \(i_1 < i_2 <
    i_3 < i_4 < i_5\). Put \(a := \lambda_{i_1}\) and \(b := \lambda_{i_2}\). The edges
    \(i_1i_3, i_1i_5 \in E(G)\) respectively imply that \(a \in B_{i_3},B_{i_5}\) (by
    Item 1 in Construction \ref{const:DS}). Similarly, the edge \(i_2i_4 \in E(G)\)
    implies that \(b \in B_{i_4}\). Observe that \(a \neq b\). Indeed, otherwise
    \(\lambda_{i_2} = b = a\), which contradicts that \(i_1i_3 \in E(G)\) (due to Item 2
    in Construction \ref{const:DS}). Now, the fact that \(a \in
    B_{i_1},B_{i_3},B_{i_5}\) and \(b \in B_{i_2},B_{i_4}\) means that the sequence is
    not \(ababa\)-free, a contradiction.

    Next, we show that \(G\) is \(\{13,35,24\}\)-free. First, we prove the following:
    \begin{claim}\label{claim:13,35,24}
        If \(i < j < k < \ell\) are vertices of \(G\) with \(ik,j\ell \in E(G)\), then
        \(f(\lambda_i) \leq f(\lambda_j)\).
    \end{claim}

    Such edges \(ik,j\ell\) are called \defn{crossing}.
    \begin{proofclaim}
        Putting \(b := \lambda_i\), \(a := \lambda_j\), suppose by contradiction that
        \(f(a) < f(b)\). In particular, \(b \neq a\). The edge \(ik\) implies that \(b
        \in B_k\), and the edge \(j\ell\) implies that \(a \in B_{\ell}\). Furthermore,
        recall that \(f(a)\) is the index of the first appearance of \(a\). Since \(f(a)
        < f(b)\), there is an appearance of \(a\) preceding the appearance of \(b =
        \lambda_i\) which was randomly chosen in \(B_i\). Together with \(a \in
        B_j,B_{\ell}\) and \(b \in B_k\), this means that the sequence is not
        \(ababa\)-free, a contradiction.
    \end{proofclaim}

    We can now show that \(G\) is \(\{13,35,24\}\)-free. Suppose by contradiction that
    \(i_1 < \dots < i_5\) form a copy of \(\{13,35,24\}\) in \(G\). By Claim
    \ref{claim:13,35,24}, we have \(f(\lambda_{i_1}) \leq f(\lambda_{i_2})\) (due to the
    crossing edges \(i_1i_3,i_2i_4\)) and \(f(\lambda_{i_2}) \leq f(\lambda_{i_3})\)
    (due to the crossing edges \(i_2i_4,i_3i_5\)). Hence, \(f(\lambda_{i_1}) \leq
    f(\lambda_{i_3})\). But since \(i_1i_3\) is an edge, Item 3 in Construction
    \ref{const:DS} gives \(f(\lambda_{i_1}) > f(\lambda_{i_3})\), a contradiction. This
    completes the proof.
\end{proof}

To complete the proof of Theorems \ref{thm:F_0} and \ref{thm:13,35,24}, it remains to
prove the following:
\begin{lem}\label{lem:DS construction independence number}
    Let \(G\) be the ordered graph given by Construction \ref{const:DS}. Then w.h.p.,
    \(\alpha(G) \leq \frac{4n}{\sqrt{h}}\).
\end{lem}

Indeed, as \(h = \Omega(\alpha(n))\), this shows that there exist \(\{13,15,24\}\)- and
\(\{13,35,24\}\)-free \(n\)-vertex ordered graphs with independence number
\(O(n/\sqrt{\alpha(n)})\), implying Theorems \ref{thm:F_0} and \ref{thm:13,35,24}.
\begin{proof}[Proof of Lemma \ref{lem:DS construction independence number}]
    The key step in the proof is to show that for every subset \(S \subseteq [n] =
    V(G)\) of size \(|S| = s\), it holds that
    \begin{equation}
        \label{eq:DS, P[S independent]}
        \mathbb{P}[S \text{ is independent}] \leq \frac{\binom{2n-1}{s}}{h^s}.
    \end{equation}
    Given \eqref{eq:DS, P[S independent]}, it is easy to complete the proof of the
    lemma. Indeed, by the union bound,
    \[
        \mathbb{P}[\alpha(G) \geq s] \leq
        \binom{n}{s} \cdot \frac{\binom{2n}{s}}{h^s} \leq
        \left( \frac{en}{s} \right)^s \cdot
        \left( \frac{2en}{s} \right)^s \cdot h^{-s} =
        \left( \frac{2e^2 n^2}{s^2 h} \right)^s = o(1)
    \]
    for \(s \geq \frac{4n}{\sqrt{h}}\).

    It remains to prove \eqref{eq:DS, P[S independent]}. So fix any \(S \subseteq [n]\)
    with \(|S| = s\). We reveal the choice of \(\lambda_i\) for \(i \in [n] \setminus
    S\) (conditioning on any given such outcome). Our goal is to upper-bound the number
    of choices of \((\lambda_i : i \in S)\) for which \(S\) is independent. To this end,
    we need the following definitions. For each symbol \(a \in [n]\), let \(S_a := \{i
    \in S : a \in B_i\}\). Namely, \(S_a\) is the set of all indices \(i \in S\) for
    which \(\lambda_i\) can equal \(a\) (as the symbol \(a\) is present in the interval
    \(B_i\)). We split the set \(S_a\) into consecutive subsets by doing the following:
    Consider all points \(j \in [n] \setminus S\) with \(\lambda_j = a\), and split
    \(S_a\) according to the intervals defined by these points \(j\). More precisely,
    let \(T_a := \{j \in [n] \setminus S : \lambda_j = a\}\) and \(t_a := |T_a|\), and
    let \(S_a^{(1)} < \dots < S_a^{(m_a)}\) be a partition of \(S_a\) such that \(m_a
    \leq t_a+1\) and for every \(1 \leq \ell \leq m_a\), there is no element of \(T_a\)
    between any two elements of \(S_a^{(\ell)}\). Note that
    \begin{equation}
        \label{eq:m_a sum}
        \sum_{a \in [n]} m_a \leq
        \sum_{a \in [n]} (t_a + 1) =
        \left( \sum_{a \in [n]} t_a \right) + n = (n - |S|) + n = 2n-s.
    \end{equation}

    For an assignment \(\lambda = (\lambda_i : i \in S)\) and for each \(a \in [n]\) and
    \(1 \leq \ell \leq m_a\), let \(c_a^{(\ell)}(\lambda)\) denote the number of \(i \in
    S_a^{(\ell)}\) with \(\lambda_i = a\). We call the sequence
    \[
        (c_a^{(\ell)}(\lambda) : a \in [n], 1 \leq \ell \leq m_a)
    \]
    the \defn{record} of \(\lambda\). The key idea is that if an assignment \(\lambda =
    (\lambda_i : i \in S)\) is such that \(S\) is independent, then \(\lambda\) is
    determined by its record.

    \begin{claim}\label{claim:DS construction, record determines assignment}
        Let \(\lambda = (\lambda_i : i \in S)\) and \(\mu = (\mu_i : i \in S)\) be two
        assignments for which \(S\) is independent. If \(c_a^{(\ell)}(\lambda) =
        c_a^{(\ell)}(\mu)\) for all \(a \in [n]\) and \(1 \leq \ell \leq m_a\), then
        \(\lambda = \mu\).
    \end{claim}
    \begin{proofclaim}
        Suppose by contradiction that \(\lambda \neq \mu\), and let \(j \in S\) be
        maximal such that \(\lambda_j \neq \mu_j\). Without loss of generality, we may
        assume that \(f(\mu_j) > f(\lambda_j)\) (else swap \(\lambda\) and \(\mu\)).
        Putting \(a := \mu_j\), let \(1 \leq \ell \leq m_a\) be such that \(j \in
        S_a^{(\ell)}\). By assumption, we have \(c_a^{(\ell)}(\lambda) =
        c_a^{(\ell)}(\mu)\). Hence, there must exist \(i \in S_a^{(\ell)}\) with
        \(\lambda_{i} = a\) and \(\mu_{i} \neq a\) (else we would have
        \(c_a^{(\ell)}(\lambda) < c_a^{(\ell)}(\mu)\)). By the maximality of \(j\), we
        have \(i \leq j\). Also, \(i \neq j\) because \(\lambda_j \neq a = \lambda_i\),
        and so \(i < j\). Now, consider a maximal \(i < j\) with \(i \in S_a^{(\ell)}\)
        and \(\lambda_i = a\). Observe that there is no \(i < k < j\) with \(\lambda_k =
        a\). Indeed, there is no such \(k \in [n] \setminus S\) because \(i,j \in
        S_a^{(\ell)}\) (and there is no element of \(T_a\) between two elements of
        \(S_a^{(\ell)}\)), and there is no such \(k \in S\) because this would force \(k
        \in S_a^{(\ell)}\) (as \(S_a^{(\ell)}\) is a consecutive subset of \(S_a\)),
        which would in turn contradict the maximality of \(i\). Finally, note that
        \(\lambda_i = a \in B_j\) because \(\mu_j = a\), and that \(f(\lambda_i) >
        f(\lambda_j)\) because \(\lambda_i = a = \mu_j\). It follows that \(ij \in
        E(G)\) for the assignment \(\lambda\), since Conditions 1-3 in Construction
        \ref{const:DS} are satisfied. This contradicts that \(S\) is independent under
        \(\lambda\).
    \end{proofclaim}

    Claim \ref{claim:DS construction, record determines assignment} implies that the
    number of choices for an assignment \(\lambda = (\lambda_i : i \in S)\) making \(S\)
    independent is at most the number of possible records of such assignments. Note that
    a record \((c_a^{(\ell)} : a \in [n], 1 \leq \ell \leq m_a)\) satisfies
    \[
        \sum_{a \in [n]} \;
        \sum_{\ell = 1}^{m_a} c_a^{(\ell)}(\lambda) = s,
    \]
    since the LHS counts each \(i \in S\) exactly once. Also, by \eqref{eq:m_a sum}, the
    number of summands on the LHS is at most \(2n-s\). Thus, the number of choices for a
    record is a balls-into-bins count with \(2n-s\) bins and \(s\) balls. The number of
    such choices is \(\binom{2n-1}{s}\). It follows that the number of assignments
    \(\lambda = (\lambda_i : i \in S)\) making \(S\) independent is at most
    \(\binom{2n-1}{s}\). On the other hand, the total number of choices for \((\lambda_i
    : i \in S)\) is \(h^s\) (as \(|B_i| = h\) for all \(i\)). This proves
    \eqref{eq:DS, P[S independent]}, completing the proof of the lemma.
\end{proof}

\section{Characterizing \(\chi\)-bounding and Ramsey-linear trees}\label{sec:charac}
We start with some preliminaries. We will need the ordered trees \(Q :=
\{13,23,24,15\}\) and \(R = \{13,23,14,25\}\) (which was defined in the introduction);
see Figure \ref{fig:Q,R}. Note that \(Q\) contains \(F_0\) and hence \(R_<(Q,K_n) =
\omega(n)\) (by Theorem \ref{thm:F_0}).
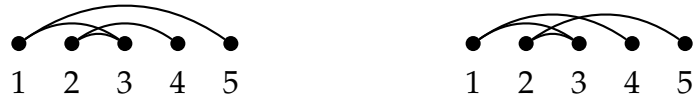
\begin{figure}[h]
    \centering
    \begin{tikzpicture}[
        vertex/.style={circle,fill=black,inner sep=2pt},
        every label/.style={label distance=3pt},
        line cap=round
        ]
        \begin{scope}
            \node[vertex,label=below:\(1\)] (c1) at (0,0) {};
            \node[vertex,label=below:\(2\)] (c2) at (0.7,0) {};
            \node[vertex,label=below:\(3\)] (c3) at (1.4,0) {};
            \node[vertex,label=below:\(4\)] (c4) at (2.1,0) {};
            \node[vertex,label=below:\(5\)] (c5) at (2.8,0) {};

            \draw[thick] (c1) to[bend left=35] (c3);
            \draw[thick] (c2) to[bend left=30] (c3);
            \draw[thick] (c2) to[bend left=35] (c4);
            \draw[thick] (c1) to[bend left=35] (c5);
        \end{scope}

        \begin{scope}[xshift=6cm]
            \node[vertex,label=below:\(1\)] (c1) at (0,0) {};
            \node[vertex,label=below:\(2\)] (c2) at (0.7,0) {};
            \node[vertex,label=below:\(3\)] (c3) at (1.4,0) {};
            \node[vertex,label=below:\(4\)] (c4) at (2.1,0) {};
            \node[vertex,label=below:\(5\)] (c5) at (2.8,0) {};

            \draw[thick] (c1) to[bend left=35] (c3);
            \draw[thick] (c2) to[bend left=30] (c3);
            \draw[thick] (c1) to[bend left=35] (c4);
            \draw[thick] (c2) to[bend left=35] (c5);
        \end{scope}
    \end{tikzpicture}
    \caption{The ordered trees \(Q\) (left) and \(R\) (right).} \label{fig:Q,R}
\end{figure}

Next, we recall two results of \cite{FKMV}.
\begin{lem}[\cite{FKMV}, Lemma 2.3]\label{lem:z-tree extremal number}
    Let \(T\) be a \(z\)-tree on \(k \geq 2\) vertices. Then \(\ex_<(n,T) = (k-2)n -
    \binom{k-1}{2}\) for all \(n \geq k\).
\end{lem}

\begin{lem}[\cite{FKMV}, Theorem 2.1]\label{lem:z-tree structural lemma}
    If an ordered tree \(T\) has interval chromatic number 2 and does not contain
    \(F_3,Q,\overleftarrow{Q},R, \overleftarrow{R}\), then \(T\) is a \(z\)-tree.
\end{lem}

In order to use Lemma \ref{lem:z-tree structural lemma}, we now show that if an ordered
tree \(T\) has no interval cut vertex and avoids certain small ordered trees, then \(T\)
has interval chromatic number 2.
\begin{lem}\label{lem:interval chromatic number 2}
    Let \(T\) be an ordered tree with no interval cut vertex, and suppose that \(T\) is
    \(F\)-free \nolinebreak for \nolinebreak all
    \begin{equation}
        \label{eq:forbidden subgraph, z-tree}
        F \in \{F_1,F_2,\overleftarrow{F_2},F_3,F_4,
        \overleftarrow{F_4},Q,\overleftarrow{Q},R,\overleftarrow{R}\}.
    \end{equation}
    Then \(\chi_<(T) = 2\).
\end{lem}
\begin{proof}
    Note that the list of forbidden subgraphs (i.e.,
    \eqref{eq:forbidden subgraph, z-tree}) is closed under order-reversal (since
    \(\overleftarrow{F}_1 = F_1\) and \(\overleftarrow{F_3} = F_3\)). Thus, we may
    always replace \(T\) with \(\overleftarrow{T}\), if necessary.

    We proceed via a sequence of claims.
    \begin{claim}\label{claim:path with endpoints to the right of all internal vertices}
        There is no path \(P = (x,v_1,\dots,v_m,y)\) in \(T\) such that \(m \geq 2\),
        and \(x,y > v_1,\dots,v_m\) or \(x,y < v_1,\dots,v_m\) in the vertex order.
    \end{claim}
    \begin{proofclaim}
        For the sake of contradiction, suppose that the claim is false, and let \(m\) be
        minimal for which there exists a counterexample. By symmetry with respect to
        reversing the order, we may assume that \(x,y > v_1,\dots,v_m\). Let \(1 \leq i
        \leq m\) be such that \(v_i\) is rightmost among \(v_1,\dots,v_m\). Observe that
        if \(i \geq 3\) then \((x,v_1,\dots,v_{i-1},v_i)\) is a smaller counterexample
        to the claim (since \(x,v_i > v_1,\dots,v_{i-1}\)); so \(i \leq 2\) by the
        minimality of \(m\). Similarly, if \(i \leq m-2\) then we get the smaller
        counterexample \((v_i,v_{i+1},\dots,v_m,y)\) (since \(y,v_i >
        v_{i+1},\dots,v_m\)), implying that \(i \geq m-1\). Hence, \(m \leq i+1 \leq
        3\).

        We now consider the cases \(m=2\) and \(m=3\) separately. Without loss of
        generality, we may assume that \(x<y\) (else traverse \(P\) in the opposite
        direction).

        Suppose first that \(m=2\), so that \(P = (x,v_1,v_2,y)\). If \(v_1 < v_2\) then
        \(v_1 < v_2 < x < y\) form a copy of \(F_4 = \{12,13,24\}\), and if \(v_2 <
        v_1\) then \(v_2 < v_1 < x < y\) form a copy of \(F_2 = \{12,23,14\}\). In
        either case, we got a contradiction to the assumption of the lemma.

        Now suppose that \(m=3\), so that \(P = (x,v_1,v_2,v_3,y)\). We saw above that
        \(i \leq 2\) and \(i \geq m-1 = 2\), and so \(i = 2\). By the choice of \(i\),
        this means that \(v_2 > v_1,v_3\). If \(v_1 < v_3\) then \(v_1 < v_3 < v_2 < x <
        y\) form a copy of \(R = \{13,23,14,25\}\), and if \(v_3 < v_1\) then \(v_3 <
        v_1 < v_2 < x < y\) form a copy of \(Q = \{ 13,23,24,15 \}\). In either case, we
        got a contradiction to the assumption of the lemma.
    \end{proofclaim}

    \begin{claim}\label{claim:no monotone 2-path}
        \(T\) contains no monotone path of length 2.
    \end{claim}
    \begin{proofclaim}
        Suppose by contradiction that \(u < v < w\) form a path (meaning that \(uv,vw
        \in E(T)\)). We claim that \(v\) is an interval cut vertex, which would
        contradict the assumption of the lemma that \(T\) has no interval cut vertex. To
        that end, it suffices to show that for every connected component \(C\) of
        \(T-v\), it holds that \(C < v\) or \(C > v\). Indeed, this would show that
        there is no edge with one endpoint to the left of \(v\) and one endpoint to the
        right of \(v\).

        So let \(C\) be any component of \(T-v\). Let \(z\) be the unique neighbor of
        \(v\) in \(C\). We show that if \(z < v\) then \(C < v\) and if \(z > v\) then
        \(C > v\). By symmetry, let us assume that \(z < v\). Suppose by contradiction
        that there is a vertex of \(C\) to the right of \(v\). Let \((z =
        z_1,z_2,\dots,z_k)\) be a path in \(C\) such that \(z_k > v\) and
        \(z_1,\dots,z_{k-1} < v\). (Namely, we traverse a path in \(C\) from \(z_1\) to
        a vertex to the right of \(v\), and stop at the first occurrence of a vertex to
        the right of \(v\).) Note that \(z_k \neq w\). Indeed, \(v\) has exactly one
        neighbor in each component of \(T-v\), and its unique neighbor in \(C\) is
        \(z\). Also, \(w\) is a neighbor of \(v\) and \(z < v < w\), implying that \(w
        \notin C\) and hence \(z_k \neq w\) (as \(z_k \in C\)). Now, consider the path
        \(P = (w,v,z=z_1,z_2,\dots,z_{k-1},z_k)\). Note that \(P\) has at least \(2\)
        internal vertices, as \(v,z\) are distinct internal vertices of \(P\). Also, the
        endpoints \(w,z_k\) of \(P\) are to the right of all internal vertices (since
        \(z_1,\dots,z_{k-1} < v < z_k,w\)). This contradicts Claim
        \ref{claim:path with endpoints to the right of all internal vertices}.

        We thus proved that if \(z < v\) then \(C < v\) and if \(z > v\) then \(C > v\),
        as required. This completes the proof of the claim.
    \end{proofclaim}
    We may now complete the proof of the lemma. Since \(T\) has no monotone path of
    length 2, for every \(v \in V(T)\) it holds that \(v\) has no left neighbors or
    \(v\) has no right neighbors. So let \(L\) (resp.~\(R\)) be the set of all \(v \in
    V(T)\) with no left neighbors (resp.~right neighbors). Then every edge of \(T\) is
    of the form \(xy\) with \(x<y\), \(x \in L\) and \(y \in R\).

    To show that \(\chi_<(T) = 2\), it now suffices to prove that \(L < R\). So suppose
    by contradiction that there exist \(x \in L\) and \(y \in R\) with \(x > y\).
    Without loss of generality, we may assume that \(x = \max(L)\) is the rightmost
    element of \(L\). Let \(P\) be the \(xy\)-path in \(T\). Let \(z\) be the rightmost
    vertex of \(P\). Since \(x \in L\), the neighbor of \(x\) on \(P\) is to the right
    of \(x\), and hence also \(z > x > y\). In particular, \(z\) is an internal vertex
    of \(P\). Let \(u,v\) be the neighbors of \(z\) on \(P\). Note that \(z \in R\),
    because \(x\) is the rightmost vertex of \(L\) and \(z > x\). Therefore, \(u,v \in
    L\), and so \(u,v \leq x\) (again, using \(x = \max(L)\)). Without loss of
    generality, suppose that \(u<v\), so that \(u < x\). Then \(u\) is an internal
    vertex of \(P\) (because \(u < x\), and \(u \in L\) while \(y \in R\)). Let \(w\) be
    the other neighbor of \(u\) on \(P\) (i.e., other than \(z\)). Note that \(u < w\)
    because \(u \in L\). Thus, we have \(u < v,w < z\) and \(uz,vz,uw \in E(T)\). Now,
    if \(v < w\) then \(u < v < w < z\) form a copy of \(F_3 = \{13,24,14\}\), and if
    \(w < v\) then \(u < w < v < z\) form a copy of \(F_1 = \{12,34,14\}\). In either
    case, we got a contradiction to the assumption of the lemma, completing the proof.
\end{proof}

By combining Lemmas \ref{lem:z-tree structural lemma} and
\ref{lem:interval chromatic number 2}, we obtain the following:
\begin{cor}\label{cor:z-tree structural lemma}
    Let \(T\) be an ordered tree with no interval cut vertex, and suppose that \(T\) is
    \(F\)-free for all \(F\) in \eqref{eq:forbidden subgraph, z-tree}. Then \(T\) is a
    \(z\)-tree.
\end{cor}

By combining Lemma \ref{lem:z-tree extremal number}, Corollary
\ref{cor:z-tree structural lemma} and the results of Section \ref{sec:prelim}, we can
easily derive Theorem \nolinebreak \ref{thm:charac}.
\begin{proof}[Proof of Theorem \ref{thm:charac}]
    The implication \(1 \Rightarrow 2\) is immediate. Let us now prove the implication
    \(3 \Rightarrow 1\). First, observe that every \(z\)-tree \(T'\) is
    \(\chi\)-bounding. Indeed, assuming \(k := |V(T')| \geq 3\), Lemma
    \ref{lem:z-tree extremal number} implies that every \(T'\)-free ordered graph has
    average degree (strictly) less than \(2k-4\). Hence, every \(T'\)-free ordered graph
    is \((2k-5)\)-degenerate and hence \((2k-4)\)-colorable. (The case \(k=2\) is
    trivial.)

    Now, let \(T\) be an ordered tree such that \(T[I]\) is a \(z\)-tree for every
    segment \(I\) of \(T\). Then by the preceding paragraph, \(T[I]\) is
    \(\chi\)-bounding for every segment \(I\) of \(T\). Now Item 1 of Corollary
    \ref{cor:segment chi-bounding} gives that \(T\) is \(\chi\)-bounding.

    Finally, we prove the implication \(2 \Rightarrow 3\). By Corollary
    \ref{cor:z-tree structural lemma}, if \(T\) is an ordered tree such that \(T[I]\) is
    not a \(z\)-tree for some segment \(I\) of \(T\), then \(T\) contains one of the
    ordered trees \(F\) from \eqref{eq:forbidden subgraph, z-tree}. For each of these
    ordered trees \(F\) we know that \(R_<(F,K_n) = \omega(n)\) (by the results of
    Section \ref{subsec:4-vertex} and by \(R_<(Q,K_n) \geq R_<(F_0,K_n) = \omega(n)\)).
    Hence, \(R_<(T,K_n) = \omega(n)\). This proves the contrapositive of the implication
    \(2 \Rightarrow 3\), as required.
\end{proof}

\section{Characterizing Ramsey-good trees}\label{sec:Ramsey-good}
Our goal is to show that an ordered tree \(T\) is Ramsey-good if and only if \(T[I]\) is
a one-sided star for every segment \(I\) of \(T\). By Item 2 in Corollary
\ref{cor:segment chi-bounding}, it suffices to consider the case that \(T\) has one
segment, i.e., that \(T\) has no interval cut vertices. By Theorem \ref{thm:charac},
\(R_<(T,K_n) = \omega(n)\) unless \(T\) is a \(z\)-tree. Clearly, if \(R_<(T,K_n) =
\omega(n)\) then \(T\) is not Ramsey-good. Thus, the remaining task is to show that the
only Ramsey-good \(z\)-trees are one-sided stars. This will be done in two steps: First
we will show that all \(z\)-trees, except for the so-called \defn{crossing
double-stars}, are not Ramsey-3-good (recall that \(T\) is Ramsey-\(n\)-good if
\(R_<(T,K_n) = (k-1)(n-1)+1\), where \(k := |V(T)|\)). We will then separately show that
if \(T\) is a crossing double-star which is not a one-sided star then \(T\) is not
Ramsey-4-good.

Recall that a \(z\)-tree \(T\) consists of an increasing tree \(T'\), say with leftmost
vertex \(\ell\) and rightmost vertex \(r\), and of additional sets of leaves \(A,B\)
such that \(B\) is adjacent to \(\ell\) and \(B > r\), and similarly \(A\) is adjacent
to \(r\) and \(A < \ell\). We say that \(T\) is a \defn{crossing double-star} if \(T'\)
is an edge (i.e., \(E(T') = \{\ell r\}\)).
\begin{defi}\label{def:crossing double star}
    A \defn{crossing double-star} is an ordered tree \(T\) with vertex-set \(\{\ell,r\}
    \cup A \cup B\), where \(A < \ell < r < B\), and with edge-set \(\{\ell r\} \cup
    \{\ell b : b \in B\} \cup \{r a : a \in A\}\).
\end{defi}

Note that if \(A\) or \(B\) is empty then \(T\) is a one-sided star. Thus, the crossing
double-stars include the one-sided stars.

We will need the following simple lemma, showing that if \(T\) is a \(z\)-tree which is
not a crossing double-star, then \(T\) or \(\overleftarrow{T}\) contains a certain
useful structure.
\begin{lem}\label{lem:z-tree not crossing double-star}
    Let \(T\) be a \(z\)-tree which is not a crossing double-star. Then in \(T\) or
    \(\overleftarrow{T}\), there exist vertices \(x < y < z\) such that
    \begin{enumerate}
        \item \(xy,xz \in E(T)\).
        \item \(x\) is not the leftmost vertex of \(T\).
        \item \(y,z\) are consecutive in the vertex order.
        \item \(x\) is adjacent to every vertex \(b \in V(T)\) with \(b > z\), and \(z\)
        is adjacent to every vertex \(a \in V(T)\) with \(a < x\).
    \end{enumerate}
\end{lem}
\begin{proof}
    We consider two cases:
    \paragraph{Case 1:} \(T\) is increasing. By definition, this means that \(T\) is
    obtained from a single edge by repeatedly adding a leftmost vertex which is a leaf
    neighbor of the rightmost vertex, or adding a rightmost vertex which is a leaf
    neighbor of the leftmost vertex. We will denote the former operation by \(L\) and
    the latter by \(R\). Observe that each of the operations \(L,R\) must be used at
    least once in the construction sequence of \(T\), because otherwise \(T\) is a
    one-sided star (which is impossible, since we assume that \(T\) is not a crossing
    double-star). By passing to \(\overleftarrow{T}\), if necessary, we may assume that
    the last operation is an \(L\). Consider the last point in the construction sequence
    of \(T\) where the operation \(R\) was used, and let \(T'\) denote the tree before
    this operation. Let \(x,y\) be the leftmost and rightmost vertices of \(T'\),
    respectively. The operation \(R\) adds a new vertex \(z > y\) which is adjacent to
    \(x\). The subsequent operations, all of which are \(L\), add vertices \(a < x\)
    which are adjacent to \(z\). Moreover, by our assumption, there is at least one such
    \(L\) operation. It is now easy to see that the conclusion of the lemma holds. (Note
    that in Case 1, \(z\) is the rightmost vertex of \(T\).)

    \paragraph{Case 2:} \(T\) is not increasing. By the definition of a \(z\)-tree,
    \(T\) consists of an increasing tree \(T'\) and leaf-sets \(A,B\) with \(A < V(T') <
    B\), such that \(A\) is adjacent to the rightmost vertex of \(T'\) and \(B\) is
    adjacent to the leftmost vertex of \(T'\). Observe that if one of \(A,B\) is empty,
    then \(T\) is increasing, contradicting the assumption of Case 2. Thus, \(A,B \neq
    \emptyset\). Also, since \(T\) is not a crossing double-star, we have \(|V(T')| \geq
    3\) (since otherwise \(T'\) would be an edge).

    Now consider the operations \(L,R\) generating the increasing tree \(T'\). By
    passing to \(\overleftarrow{T}\), if necessary, we may assume that the last
    operation in the construction of \(T'\) is \(R\). Thus, this operation adds a
    rightmost leaf \(z\) (of \(T'\)), adjacent to the leftmost vertex \(x\) of \(T'\).
    Let \(y\) be the rightmost vertex of \(T' - \{z\}\). Then \(x,y\) are adjacent,
    because the leftmost and rightmost vertices of any increasing tree are adjacent
    (here we use this for the tree \(T' - \{z\}\)). It is now easy to see that the
    conclusion of the lemma holds. Indeed, \(x\) is not the leftmost vertex of \(T\)
    because \(A \neq \emptyset\); and \(x\) is adjacent to all vertices in \(B\) and
    \(z\) is adjacent to all vertices in \(A\).
\end{proof}

\begin{lem}\label{lem:not Ramsey-3-good}
    Let \(T\) be a \(k\)-vertex \(z\)-tree which is not a crossing double-star. Then
    \(R_<(T,K_3) > 2k-1\) (and hence \(T\) is not Ramsey-3-good).
\end{lem}
\begin{proof}
    By Lemma \ref{lem:z-tree not crossing double-star}, by passing to
    \(\overleftarrow{T}\) if necessary, we may write \(V(T)\) as
    \[
        A < x < C < y < z < B
    \]
    such that \(A \neq \emptyset\), and the following are edges:
    \[
        xy, xz,
        \qquad
        \{x b : b \in B\},
        \qquad
        \{z a : a \in A\}.
    \]

    To prove the lemma, we need to construct an ordered graph on \(2k-1\) vertices which
    is \(T\)-free and has no independent set of size 3. The construction will be a
    clique-blowup of \(C_5\). Namely, it will consist of 5 vertex-sets \(V_1,\dots,V_5\)
    such that each \(V_i\) is a clique and \((V_i,V_{i+1})\) is a complete bipartite
    graph for every \(1 \leq i \leq 5\), with indices taken modulo \(5\). It is easy to
    see that any such graph has independence number 2 (as \(\alpha(C_5) = 2\) and the
    sets \(V_1,\dots,V_5\) are cliques). It remains to choose the sizes of the sets
    \(V_1,\dots,V_5\) and the vertex order so as to avoid \(T\) and have \(2k-1\)
    vertices in \nolinebreak total. \nolinebreak Set
    \[
        |V_1| = |V_4| = 1,
        \qquad
        |V_2| = k-|B|-2,
        \qquad
        |V_3| = |B|+1.
    \]
    Also, split \(V_5\) into two sets \(U_5,W_5\) with
    \[
        |U_5| = k-|B|-2,
        \qquad |W_5| = |B|.
    \]
    Then
    \[
        |V_1| + \dots + |V_5| = 1 + (k-|B|-2) + (|B|+1) + 1 + (k-2) = 2k-1,
    \]
    as required.

    Next, order the vertices as follows:
    \[
        V_1,V_2,U_5,V_3,V_4,W_5.
    \]
    We claim that the resulting ordered graph \(G\) is \(T\)-free. So suppose by
    contradiction that there is an embedding \(\varphi : T \rightarrow G\). We split
    into cases according to the location of \(\varphi(x)\). First, observe that
    \(\varphi(x) \in V_1\) is impossible, because \(|V_1| = 1\) and \(\varphi(A) <
    \varphi(x)\), meaning that \(\varphi(x) \in V_1\) leaves no suitable image for the
    vertices of \(A\) (here we use that \(A \neq \varnothing\)). Also, note that \(x\)
    has at least \(|B|+2\) right neighbors in \(T\) (namely, the vertices in \(\{y,z\}
    \cup B\)), whereas every vertex in \(V_3 \cup V_4 \cup W_5\) has at most \(|B|+1\)
    right neighbors in \(G\). It follows that \(\varphi(x) \notin V_3,V_4,W_5\).

    The above only leaves the cases \(\varphi(x) \in V_2\) and \(\varphi(x) \in U_5\).

    \paragraph{Case 1:} \(\varphi(x) \in V_2\). Then all right neighbors of
    \(\varphi(x)\) are in \(V_2 \cup V_3\). Also, \(|V_3| = |B|+1\), and \(\varphi(x)\)
    has \(|B|+1\) neighbors to the right of \(\varphi(y)\) (namely, the images of
    \(z,B\)). Hence, \(\varphi(y)\) cannot be in \(V_3\), and so \(\varphi(y) \in V_2\).
    It follows that the entire set \(\{x,y\} \cup A \cup C\) is mapped to \(V_1 \cup
    V_2\) (since \(x,A,C < y\)). Now, note that \(2 + |A| + |C| = |V(T)| - 1 - |B| =
    k-|B|-1 = |V_1 \cup V_2|\). Thus, the set \(V_1 \cup V_2\) is precisely the image of
    \(\{x,y\} \cup A \cup C\). In particular, this means that \(\varphi(z) \in V_3\),
    and that the first vertex of \(A\) is mapped to \(V_1\). But \(z\) is adjacent to
    all vertices of \(A\), while \(V_3\) is not adjacent to \(V_1\), a contradiction.

    \paragraph{Case 2:} \(\varphi(x) \in U_5\). Then all right neighbors of
    \(\varphi(x)\) are in \(V_4 \cup V_5\). We claim that \(\varphi(y) \in U_5\).
    Indeed, otherwise, \(\varphi(y) \in V_4 \cup W_5\). But then all of \(y,z,B\) are
    mapped to \(V_4 \cup W_5\), whereas \(|V_4 \cup W_5| = |B|+1 < |B|+2\), a
    contradiction. Thus, \(\varphi(y) \in U_5\), and hence also \(\varphi(C) \subseteq
    U_5\) (since \(x < C < y\)).

    Note that \(\varphi(z) \in V_4 \cup V_5\). Furthermore, \(\varphi(z) \in W_5\) is
    impossible, since \(\varphi(x)\) needs to have \(|B|\) neighbors to the right of
    \(\varphi(z)\), whereas \(|W_5| = |B|\). Thus, \(\varphi(z) \in U_5 \cup V_4\).
    Moreover, each vertex in \(\varphi(A)\) is adjacent to \(\varphi(z)\) and is to the
    left of \(\varphi(x)\). Therefore, \(\varphi(A) \subseteq V_1 \cup U_5\).

    We now consider two cases: If \(\varphi(z) \in U_5\), then the entire set
    \(\{x,y,z\} \cup A \cup C\), which has size \(k-|B|\), is mapped into \(V_1 \cup
    U_5\), which has size \(k-|B|-1\), a contradiction. And if \(\varphi(z) \in V_4\),
    then \(\varphi(A) \cap V_1 = \varnothing\) (since the vertices of \(\varphi(A)\)
    need to be adjacent to \(\varphi(z)\)), meaning that \(\{x,y\} \cup A \cup C\) is
    mapped into \(U_5\). This is similarly impossible, since \(\{x,y\} \cup A \cup C\)
    has size \(k-|B|-1\), while \(|U_5| = k-|B|-2\). This completes \nolinebreak the
    \nolinebreak proof.
\end{proof}

It remains to handle the crossing double-stars. The following lemma states that if \(T\)
is a crossing double-star but not a one-sided star then \(T\) is not Ramsey-4-good.
\begin{lem}\label{lem:not Ramsey-4-good}
    Let \(T\) be a \(k\)-vertex crossing double-star which is not a one-sided star. Then
    \(R_<(T,K_4) > 3k-2\).
\end{lem}
\begin{proof}
    By definition, \(V(T)\) can be written as \(A < \ell < r < B\), where \(\ell r \in
    E(T)\); \(\ell\) is adjacent to all vertices of \(B\); and \(r\) is adjacent to all
    vertices of \(A\). We have \(A,B \neq \emptyset\), because otherwise \(T\) is a
    one-sided star. Put \(a := |A|, b := |B|\), so that \(a+b = k-2\).

    Our goal is to construct a graph on \(3k-2\) vertices which is \(T\)-free and has no
    independent set of size \(4\). The construction is a clique-blowup of the following
    graph \(G_0\) on the vertices \(1,\dots,8\):
    \[
        E(G_0) = \{ 12,23,34,45,56,67,78,81 \} \cup \{ 17,48 \}.
    \]
    Thus, \(G_0\) is obtained from the 8-cycle \((1,2,\dots,8,1)\) by adding the chords
    \(17,48\). It is straightforward to check that \(\alpha(G_0) = 3\) (indeed, the
    \(8\)-cycle has exactly two independent sets of size 4, and both are blocked by the
    chords \(17,48\)). Hence, each clique-blowup of \(G_0\) also has independence number
    3.

    Now consider the following clique-blowup of \(G_0\): The blowup sets are
    \(V_1,\dots,V_8\), and
    \[
        |V_1| = |V_3| = |V_4| = |V_6| = 1,
        \qquad
        |V_7| = b,
        \qquad
        |V_8| = a.
    \]
    Also, split \(V_2\) and \(V_5\) into sets \(U_2,W_2\) and \(U_5,W_5\), respectively,
    and put
    \[
        |U_2| = |U_5| = a,
        \qquad
        |W_2| = |W_5| = b.
    \]
    The total number of vertices of the resulting clique-blowup \(G\) is \( 4 + a + b +
    2(a+b) = 3(a+b) + 4 = 3k-2, \) as required.

    The vertex ordering of \(G\) is
    \[
        U_2,V_1,V_8,V_4,U_5,V_7,W_2,V_6,W_5,V_3.
    \]
    To prove that \(G\) is \(T\)-free, we need to show that for every edge \(uv \in
    E(G)\) with \(u<v\), either \(v\) has fewer than \(a\) neighbors to the left of
    \(u\), or \(u\) has fewer than \(b\) neighbors to the right of \(v\). Indeed, this
    would exactly mean that there is no embedding \(\varphi : T \rightarrow G\) with
    \(\varphi(\ell)=u, \varphi(r) = v\). Proving this statement is a somewhat tedious
    case analysis.

    First, consider the union \(L = U_2 \cup V_1 \cup V_8 \cup V_4 \cup U_5\) of the 5
    left blocks in the vertex ordering. The pairs of adjacent blocks among these 5 are
    \((U_2,V_1),(V_1,V_8),(V_8,V_4),(V_4,U_5)\). It is now easy to check that due to the
    sizes of the blocks (i.e., \((|U_2|,|V_1|,|V_8|,|V_4|,|U_5|) = (a,1,a,1,a)\)), every
    vertex \(x \in L\) has at most \(a\) left neighbors. This means that if \(u,v \in
    L\) with \(u < v\) and \(uv \in E(G)\) then \(v\) has at most \(a-1\) neighbors to
    the left of \(u\). Thus, this proves the required claim in the case \(u,v \in L\).

    Similarly, let \(R = V_7 \cup W_2 \cup V_6 \cup W_5 \cup V_3\) be the union of the 5
    right blocks. The pairs of adjacent blocks among these 5 are
    \((V_7,V_6),(W_2,V_3),(V_6,W_5)\). It is now easy to check that due to the sizes of
    the blocks (i.e., \((|V_7|,|W_2|,|V_6|,|W_5|,|V_3|) = (b,b,1,b,1)\)), every vertex
    \(x \in R\) has at most \(b\) right neighbors. This means that if \(u,v \in R\) with
    \(u < v\) and \(uv \in E(G)\), then \(u\) has at most \(b-1\) neighbors to the right
    of \(v\), as required.

    It remains to rule out the case that \(u \in L, v \in R\). In that case, \(u \in X,
    v \in Y\), where \((X,Y)\) equals one of the pairs
    \[
        (U_2,W_2), (U_2,V_3), (V_1,V_7), (V_1,W_2), (V_8,V_7), (V_4,W_5), (V_4,V_3),
            (U_5,V_6),
        (U_5,W_5).
    \]
    It is easy to see that for each of these pairs \((X,Y)\), either \(X\) is the
    leftmost block in \(L\) adjacent to \(Y\), or \(Y\) is the rightmost block in \(R\)
    adjacent to \(X\). Now, since \(|X| \leq a\) and \(|Y| \leq b\) (by the choice of
    the block-sizes), either \(v\) has at most \(a-1\) neighbors left of \(u\), or \(u\)
    has at most \(b-1\) neighbors right of \(v\), as required.
\end{proof}

Given all of the above, Theorem \ref{thm:Ramsey-good} now easily follows:
\begin{proof}[Proof of Theorem \ref{thm:Ramsey-good}]
    By Item 2 of Corollary \ref{cor:segment chi-bounding}, it suffices to show that for
    every tree \(T\) with no interval cut vertex, it holds that \(T\) is Ramsey-good if
    and only if \(T\) is a one-sided star. First, every one-sided star is Ramsey-good.
    Indeed, let \(n \geq 2\), and let \(G\) be an ordered graph on \((k-1)(n-1)+1\)
    vertices with no independent set of size \(n\). Then there must exist a vertex with
    forward degree at least \(k-1\), since otherwise an independent set of size \(n\)
    can be constructed greedily by repeatedly adding the leftmost vertex to the
    independent set and removing its forward neighbors. By symmetry, there also exists a
    vertex with backward degree at least \(k-1\). This shows that \(G\) contains both
    one-sided stars on \(k\) vertices.

    Conversely, let \(T\) be an ordered tree with no interval cut vertex, and suppose
    that \(T\) is Ramsey-good. Then in particular \(R_<(T,K_n) = O(n)\), and so \(T\) is
    a \(z\)-tree by Theorem \ref{thm:charac}. By Lemmas \ref{lem:not Ramsey-3-good} and
    \ref{lem:not Ramsey-4-good}, \(T\) is a one-sided star.
\end{proof}

\section{Additional positive results}\label{sec:reductions}

\subsection{Proof of Theorem \ref{thm:insertion-into-stars}}

We first require the following auxiliary lemma.

\begin{lem}\label{lem:k-shift graph}
    Let \(G\) be a graph and let \(V(G) = U_1 \cup \dots \cup U_t\) be a partition into
    independent sets. Let \(k \geq 0\), and suppose that there is an orientation of
    \(E(G)\) such that for every \(x \in V(G)\), there are at most \(k\) indices \(i \in
    [t]\) such that \(x\) has out-neighbors in \(U_i\). Then \(\alpha(G) \geq
    \frac{1}{e(k+1)}|V(G)|\).
\end{lem}
\begin{proof}
    When \(k = 0\), \(G\) is edgeless and the result is immediate. Hence, assume that
    \(k \ge 1\). For each \(x \in V(G)\), let \(I(x)\) be the set of all \(i \in [t]\)
    such that \(x\) has out-neighbors in \(U_i\) in the given orientation. By
    assumption, \(|I(x)| \leq k\). Pick a subset \(J \subseteq [t]\) randomly by
    including each \(i \in [t]\) in \(J\) independently with probability \(p :=
    \frac{1}{k+1}\). Let \(S\) be the set of all \(x \in V(G)\) such that, letting \(i
    \in [t]\) be such that \(x \in U_i\), it holds that \(i \in J\) and \(I(x) \cap J =
    \emptyset\). Note that \(S\) is independent for every choice of \(J\). Indeed,
    suppose by contradiction that \(x,y \in S\) for some edge \(xy \in E(G)\). Without
    loss of generality, suppose that this edge is oriented from \(x\) to \(y\). Let
    \(i,j \in [t]\) be such that \(x \in U_i, y \in U_j\). By definition, \(j \in
    I(x)\). Now, \(y \in S\) implies that \(j \in J\), whereas \(x \in S\) implies that
    \(j \notin J\), a contradiction.

    For every \(x \in V(G)\), we have
    \[
        \mathbb{P}[x \in S] = p(1-p)^{|I(x)|} \geq
        p(1-p)^k =
        \frac{1}{k+1} \cdot \left( 1 - \frac{1}{k+1} \right)^k \geq \frac{1}{e(k+1)}.
    \]
    Hence, by linearity of expectation, \(\mathbb{E}[|S|] \geq \frac{1}{e(k+1)}|V(G)|\).
    This completes the proof.
\end{proof}

\begin{rem}
    In the setting of Lemma \ref{lem:k-shift graph}, if \(G\) is an ordered graph,
    \(U_1,\dots,U_t\) are intervals and \(k=1\), then the resulting structure is closely
    related to the so-called shift graphs. The shift graph \(S_n\) has vertices \(V(S_n)
    = \binom{[n]}{2}\), ordered lexicographically, and \(\{(x,y),(y,z)\}\) is an edge
    for every triple \(1 \leq x < y < z \leq n\). Taking \(U_x := \{(x,y) : x < y \leq
    n\}\) for \(1 \leq x \leq n-1\) and orienting all edges forward, it is easy to see
    that every vertex has out-neighbors in at most one interval among
    \(U_1,\dots,U_{n-1}\). Conversely, one can show that every ordered graph having such
    structure is contained in the blowup of some shift graph, possibly after rearranging
    the vertices in the independent sets \(U_i\).
\end{rem}

\begin{proof}[Proof of Theorem \ref{thm:insertion-into-stars}]
    Recall that \(H_1,\dots,H_s\) are ordered graphs with \(R_<(H_i,K_n) = O(n)\) for
    all \(i \in [s]\), and \(F\) is the ordered graph with vertex-set
    \[
        x < V(H_1) < y_1 < V(H_2) < y_2 < \dots < V(H_s) < y_s,
    \]
    consisting of the graphs \(H_1,\dots,H_s\) and of the star with center \(x\) and
    leaves \(y_1,\dots,y_s\). Our goal is to show that \(R_<(F,K_n) = O(n)\).

    Let \(H\) be the ordered graph consisting of disjoint copies of \(H_1,\dots,H_s\)
    with \(V(H_1) < \dots < V(H_s)\). We may assume that \(|V(H)| \geq 2\), because
    otherwise \(F\) is contained in a one-sided star, and hence \(R_<(F,K_n) = O(n)\) by
    Theorem \ref{thm:Ramsey-good}. Since \(R_<(H,K_n) \leq \sum_{i=1}^s R_<(H_i,K_n) =
    O(n)\), there is a constant \(c \geq 1\) such that every \(H\)-free ordered graph on
    \(m\) vertices has an independent set of size at least \(m/c\).

    Let \(G\) be an \(F\)-free ordered graph on \(N\) vertices. Partition \(V(G)\) into
    intervals \(I_1 < \dots < I_t\), where \(I_i\) is the longest initial segment of
    \(V(G) \setminus (I_1 \cup \dots \cup I_{i-1})\) inducing an \(H\)-free graph (we
    can do this because a 1-vertex interval is \(H\)-free, as \(|V(H)| \geq 2\)). By
    maximality, \(G[I_i \cup I_{i+1}]\) contains a copy of \(H\) for every \(i < t\).
    Hence, if \(p \in I_i\) and \(q \in I_j\) with \(j \geq i+3\), then \(I_{i+1} \cup
    I_{i+2}\) contains a copy of \(H\), which lies strictly between \(p\) and \(q\).

    We now claim that every vertex \(x\) has forward neighbors in at most \(3s\) of the
    intervals \(I_1,\dots,I_t\). Otherwise, there are indices \(j_0 < j_1 < \dots <
    j_{3s}\) such that \(x \in I_{j_0}\) and \(x\) has a forward neighbor in each of
    \(I_{j_1},\dots,I_{j_{3s}}\). Let \(y_0 := x\) and, for each \(r \in [s]\), let
    \(y_r\) be a forward neighbor of \(x\) in \(I_{j_{3r}}\). Since \(j_{3r} \geq
    j_{3(r-1)}+3\), there is a copy of \(H\), and hence of \(H_r\), strictly between
    \(y_{r-1}\) and \(y_r\) (for each \(1 \leq r \leq s\)). These copies together with
    \(x,y_1,\dots,y_s\) form a copy of \(F\), a contradiction.

    Finally, for each \(i \in [t]\), fix an independent set \(U_i \subseteq I_i\) of
    size at least \(|I_i|/c\) (such an independent set exists because \(G[I_i]\) is
    \(H\)-free). Applying Lemma \ref{lem:k-shift graph} to \(G[U_1 \cup \dots \cup
    U_t]\) with the partition \(U_1,\dots,U_t\), all edges oriented forward, and \(k =
    3s\), we obtain \(\alpha(G) \geq \frac{N}{ec(3s+1)}\). Thus, \(R_<(F,K_n) = O(n)\).
\end{proof}

\subsection{Proof of Theorem \ref{thm:one-vertex interlace}}

\begin{proof}[Proof of Theorem \ref{thm:one-vertex interlace}]
    Write \(V(A) = \{a_1 < \dots < a_k\}\) and \(V(B) = \{b_1 < \dots < b_{\ell}\}\), so
    that \(A \diamond B\) has vertices \(a_1 < \dots < a_{k-1} < b_1 < a_k < b_2 < \dots
    < b_{\ell}\). Let \(G\) be a graph with no copy of \(A \diamond B\) and with
    \(\alpha(G) < n\), and let us show that
    \begin{equation}
        \label{eq:one-vertex interlace}
        |V(G)| \leq 2R_<(A,K_n) + 2R_<(B,K_n)-4.
    \end{equation}

    Setting \(N := |V(G)|\), identify \(V(G)\) with \([N]\). Let \(X\) be the set of all
    vertices \(x \in V(G)\) such that there exists {\bf no} copy of \(B\) with vertices
    \(x = x_1 < x_2 < \dots < x_{\ell}\) (with \(x_i\) playing the role of \(b_i\) for
    \(i=1,\dots,\ell\)) and such that \(x_2 \geq x_1+2\). For \(j=0,1\), let \(X_j\)
    denote the set of all \(x \in X\) with \(x \equiv j \pmod{2}\). Observe that
    \(G[X_j]\) is \(B\)-free. Indeed, a copy \(x_1 < \dots < x_{\ell}\) of \(B\) with
    \(x_1,\dots,x_{\ell} \in X_j\) clearly has \(x_2 \geq x_1+2\), which would imply
    that \(x_1 \notin X\), a contradiction. Now, since \(\alpha(G) < n\), we conclude
    that \(|X_j| \leq R_<(B,K_n)-1\), and hence
    \begin{equation}
        \label{eq:one-vertex interlace, aux 1}
        |X| \leq 2R_<(B,K_n)-2.
    \end{equation}

    Next, let \(X^+ := \{x+1 : x \in X \setminus \{N\}\}\); namely, \(X^+\) is the set
    of immediate successors of the vertices in \(X\). Also, set \(Y := V(G) \setminus
    X^+\) and, for \(j = 1,2\), let \(Y_j\) be the set of all \(y \in Y\) with \(y
    \equiv j \pmod{2}\). We claim that \(G[Y_j]\) is \(A\)-free. Indeed, suppose by
    contradiction that \(y_1 < \dots < y_k\) are vertices of \(Y_j\) forming a copy of
    \(A\) (with \(y_i\) playing the role of \(a_i\) for \(i=1,\dots,k\)). Since \(y_k
    \in Y = V(G) \setminus X^+\), this means that \(x := y_k-1 \notin X\). Note also
    that \(y_{k-1} < x\) because \(y_{k-1},y_k\) have the same parity. Now, by the
    definition of \(X\) (and as \(x \notin X\)), there exists a copy \(x = x_1 < x_2 <
    \dots < x_{\ell}\) of \(B\) with \(x_2 \geq x_1 + 2 > y_k\). In particular, \(y_1 <
    \dots < y_{k-1} < x_1 < y_k < x_2 < \dots < x_{\ell}\) form a copy of \(A \diamond
    B\), a contradiction. This shows that \(Y_j\) is \(A\)-free. It follows that \(|Y| =
    |Y_1| + |Y_2| \leq 2R_<(A,K_n) - 2\). Finally, \[|V(G)| = |Y| + |X^+| \leq |Y| + |X|
    \leq 2R_<(A, K_n) + 2R_<(B, K_n) - 4.\qedhere\]
\end{proof}

\subsection{Proof of Theorem \ref{thm:chi-bounding reduction}}

\begin{proof}[Proof of Theorem \ref{thm:chi-bounding reduction}]
    Let \(F\) be an ordered graph with vertices \(v_1 < \dots < v_k\) and \(v_1v_k \in
    E(F)\). By assumption, every \(F\)-free ordered graph has chromatic number at most
    \(C\). Recall that \(F^\star\) is the ordered graph with vertices \(a < v_1 < \dots
    < v_k < b < c\) consisting of the edges of \(F\) and the edges \(ab,v_kc\). Let
    \(G\) be an \(F^\star\)-free ordered graph. Our goal is to show that \(\chi(G) \leq
    5C+4\).

    First, let \(X\) be the set of vertices \(x\) such that there exists a copy of \(F\)
    in \(G\) in which \(x\) plays the role of \(v_k\). Note that \(G - X := G[V(G)
    \setminus X]\) is \(F\)-free and hence \(\chi(G-X) \leq C\). It thus suffices to
    show that \(\chi(G[X]) \leq 4C+4\).

    We will need the following definitions: For each \(x \in X\), let \(f(x)\) denote
    the leftmost neighbor of \(x\) in the graph \(G\), and let \(g(x)\) denote the
    rightmost vertex \(y \in V(G)\) such that there exists a copy of \(F\) in which
    \(y\) plays the role of \(v_1\) and \(x\) plays the role of \(v_k\). Note that such
    a vertex \(y\) is adjacent to \(x\), because \(v_1v_k \in E(F)\). Since \(x \in X\),
    it follows that \(f(x),g(x)\) are well-defined and \(f(x) \leq g(x) < x\) in the
    vertex order.

    Next, let \(G_1\) be the graph on \(X\) obtained from \(G[X]\) by deleting, for each
    \(x \in X\), the edge to the rightmost forward neighbor of \(x\) (if \(x\) has
    forward neighbors). Let \(G_2\) denote the graph of deleted edges, so that \(G[X] =
    G_1 \cup G_2\) (as an edge-union). Note that \(G_2\) is 1-degenerate (i.e., a
    forest), because each \(x \in X\) has at most one forward neighbor in \(G_2\). Thus,
    \(G_2\) is 2-colorable. Therefore, in order to prove that \(\chi(G[X]) \leq 4C+4\),
    it suffices to prove that \(\chi(G_1) \leq 2C+2\).

    The key observation we need is as follows:
    \begin{claim}\label{claim:chi-bounding reduction}
        Let \(x < x' \leq x''\) be vertices of \(X\) with \(xx'' \in E(G_1)\). Then
        \(g(x) \leq f(x')\).
    \end{claim}
    \begin{proofclaim}
        Suppose by contradiction that \(y' := f(x') < g(x) =: y\). By the definition of
        \(f(x')\), we have \(y'x' \in E(G)\). By the definition of \(g(x)\), there
        exists a copy of \(F\) with vertices \(z_1 < \dots < z_k\) such that \(z_1 = y\)
        and \(z_k = x\). Finally, as \(xx'' \in E(G_1)\), the vertex \(x''\) is not the
        rightmost neighbor of \(x\). Thus, there exists \(w > x'' \geq x'\) with \(xw
        \in E(G)\). But now \(y' < z_1 < \dots < z_k < x' < w\) form a copy of
        \(F^\star\) in which \(z_1,\dots,z_k\) play the roles of \(v_1,\dots,v_k\),
        respectively; and \(y',x',w\) play the roles of \(a,b,c\), respectively. But
        this contradicts the assumption that \(G\) is \(F^\star\)-free.
    \end{proofclaim}

    In order to prove that \(\chi(G_1) \leq 2C+2\), we define a certain decomposition of
    \(V(G_1) = X\), as follows. For a non-empty subset \(Y \subseteq X\), let \(L(Y)\)
    denote the set of all \(x \in Y\) such that \(f(x) < Y\), meaning that \(x\) has a
    neighbor preceding \(Y\). Observe that if \(x_0 = \min(Y)\) denotes the leftmost
    element of \(Y\), then \(f(x_0) < Y\) because \(f(x_0) < x_0\), and so \(L(Y) \neq
    \emptyset\). In our decomposition of \(X\), we will repeatedly decompose a set \(Y\)
    into the set \(L(Y)\) and the \defn{maximal intervals of \(Y \setminus L(Y)\)}; by
    the latter, we mean the maximal subsets \(Y_1 < \dots < Y_m\) of \(Y \setminus
    L(Y)\) such that there is no element of \(L(Y)\) between any two elements of \(Y_i\)
    (for each \(1 \leq i \leq m\)), and there is an element of \(L(Y)\) between \(Y_i\)
    and \(Y_{i+1}\) (for each \(1 \leq i \leq m-1\)). These definitions have the
    following properties:
    \begin{claim}\label{claim:L(Y) edges}
        Let \(Y \subseteq X\) be a non-empty subset, and let \(Y_1 < \dots < Y_m\) be
        the maximal intervals of \(Y \setminus L(Y)\). Then:
        \begin{enumerate}
            \item There are no edges in \(G_1\) between \(Y_1,\dots,Y_m\).
            \item For every \(i \in [m]\) and every edge \(xy \in E(G_1)\) with \(x \in
            L(Y)\) and \(y \in Y_i\), it holds that \(y \in L(Y_i)\).
        \end{enumerate}
    \end{claim}
    \begin{proofclaim}
        We start with Item 1. Suppose by contradiction that there is an edge \(xy\) with
        \(x \in Y_i\), \(y \in Y_j\), and \(i < j\) (and so \(x < y\)). Let \(x' \in
        L(Y)\) be an element of \(L(Y)\) between \(Y_i\) and \(Y_j\). Note that \(f(x')
        < f(x)\), because \(x'\) has a neighbor to the left of \(Y\) while \(x\) does
        not (as \(x' \in L(Y)\), \(x \notin L(Y)\)). Also, \(f(x) \leq g(x)\) (this
        holds for every vertex). Thus, \(f(x') < g(x)\), contradicting Claim
        \ref{claim:chi-bounding reduction}.

        For Item 2, let \(xy \in E(G_1)\) with \(x \in L(Y)\) and \(y \in Y_i\). First,
        we claim that \(x < y\). Indeed, we have \(f(x) < f(y) \leq g(y)\) (as above),
        so \(y < x\) would contradict Claim \ref{claim:chi-bounding reduction}. Now,
        since \(x < y \in Y_i\) and \(Y_i\) is a maximal interval of \(Y \setminus
        L(Y)\), it follows that \(x < Y_i\). Hence, \(y\) has a neighbor preceding
        \(Y_i\) (namely \(x\)), which means, by definition, that \(y \in L(Y_i)\).
    \end{proofclaim}

    We now define the decomposition of \(X\) alluded to above. Starting with
    \(\mathcal{P}_0 = \{X\}\), do the following for each \(i \geq 1\): Initialize
    \(\mathcal{P}_i := \emptyset\). For every \(Y \in \mathcal{P}_{i-1}\), remove the
    vertices of \(L(Y)\), let \(Y_1,\dots,Y_m\) be the maximal intervals of \(Y
    \setminus L(Y)\), and add \(Y_1,\dots,Y_m\) to \(\mathcal{P}_i\). The sets
    \(Y_1,\dots,Y_m\) are called the \defn{children} of \(Y\). This defines
    \(\mathcal{P}_i\). We observed above that \(L(Y) \neq \emptyset\), and so at least
    one vertex of \(Y\) is removed for each \(Y \in \mathcal{P}_{i-1}\). This means that
    this process eventually removes all vertices of \(X\). Using Item 1 of Claim
    \ref{claim:L(Y) edges}, it is easy to show by induction that for every \(i\), there
    are no edges of \(G_1\) between any two sets in \(\mathcal{P}_i\). By Item 2 of
    Claim \ref{claim:L(Y) edges}, it follows that for every edge \(e\) of \(G_1\), there
    is an \(i \geq 0\) and a set \(Y \in \mathcal{P}_i\) such that \(e\) is either
    contained in \(L(Y)\) or goes between \(L(Y)\) and \(L(Y')\) for some child \(Y'\)
    of \(Y\). Thus, to show that \(\chi(G_1) \leq 2C+2\), it suffices to show that for
    every set \(Y\), the subgraph of \(G_1\) induced by \(L(Y)\) is \((C+1)\)-colorable.
    Indeed, we can then use one set of \(C+1\) colors to color \(L(Y)\) for all \(Y \in
    \mathcal{P}_i\) with \(i\) even, and another set of \(C+1\) colors to color \(L(Y)\)
    for all \(Y \in \mathcal{P}_i\) with \(i\) odd. Thus, it remains to prove the
    following:
    \begin{claim}\label{claim:L(Y) coloring}
        For every set \(Y \subseteq X\), the subgraph of \(G_1\) induced by \(L(Y)\) is
        \((C+1)\)-colorable.
    \end{claim}
    \begin{proofclaim}
        Let \(Z\) be the set of all \(y \in L(Y)\) with \(g(y) < Y\). First, we claim
        that \(G_1[Z]\) is \(F\)-free. Indeed, suppose by contradiction that \(G_1[Z]\)
        contains a copy of \(F\) with vertices \(z_1 < \dots < z_k\) and with \(z_i\)
        playing the role of \(v_i\) (for \(i=1,\dots,k\)). By the definition of the
        function \(g\), we then have \(g(z_k) \geq z_1\), contradicting that \(g(z_k) <
        Y\) (since \(z_1 \in L(Y) \subset Y\)). This proves that \(G_1[Z]\) is
        \(F\)-free, and hence \(\chi(G_1[Z]) \leq C\). Finally, we claim that \(L(Y)
        \setminus Z\) is independent in \(G_1\). Indeed, suppose by contradiction that
        there exists \(xy \in E(G_1)\) with \(x,y \in L(Y) \setminus Z\) and \(x < y\).
        Then \(f(y) < Y\) (because \(y \in L(Y)\)) and \(g(x) \geq \min(Y)\) (since \(x
        \notin Z\)). Hence, \(f(y) < g(x)\), contradicting Claim
        \ref{claim:chi-bounding reduction}.
    \end{proofclaim}
\end{proof}

\section{Concluding remarks and open problems}\label{sec:concluding remarks}

\begin{itemize}
    \item The problem of characterizing the ordered forests with linear off-diagonal
    Ramsey numbers remains open (see Problem~\ref{prob:Ramsey-linear}), as does the
    problem of characterizing the \(\chi\)-bounding ordered forests:
    \begin{prob}\label{prob:12,14,35}
        Characterize the \(\chi\)-bounding ordered forests.
    \end{prob}
    \item Theorem~\ref{thm:F_0} shows that for some ordered forests \(F\), the Ramsey
    number \(R_<(F,K_n)\) can be only barely superlinear. We wonder whether the same can
    hold for a tree.
    \begin{prob}\label{prob:trees superlinear}
        Is there an ordered tree \(T\) for which \(\omega(n) \leq R_<(T,K_n) \leq
        n\log^{o(1)}(n)\)?
    \end{prob}
    We note that for the extremal number \(\ex_<(n,T)\) of ordered trees, the answer is
    negative: If \(\ex_<(n,T) = \omega(n)\) then \(\ex_<(n,T) = \Omega(n\log n)\)
    (see~\cite{FKMV}).

    We suspect that \(Q = \{13,23,24,15\}\) could give an affirmative answer to
    Problem~\ref{prob:trees superlinear}. We also wonder whether a polylogarithmic gap
    from linearity necessarily implies an almost logarithmic gap:
    \begin{prob}\label{prob:trees logarithmic gap}
        Is it true that every ordered tree \(T\) satisfying \(R_<(T,K_n) \geq
        n\log^{\Omega(1)}(n)\) also satisfies \(R_<(T,K_n) \geq n\log^{1-o(1)}(n)\)?
    \end{prob}
    The results of Section \ref{subsec:4-vertex} provide some supporting evidence for a
    positive answer to Problem \ref{prob:trees logarithmic gap}. Finally, we ask which
    ordered trees satisfy this stronger lower bound.
    \begin{prob}
        Characterize the ordered trees \(T\) for which \(R_<(T,K_n) \geq
        n\log^{1-o(1)}(n)\).
    \end{prob}
    \item The precise asymptotics of \(R_<(F,K_n)\) remain undetermined for some small
    ordered forests \(F\), including \(R\) (see Theorem~\ref{thm:R}), \(F_0\) (see
    Theorem~\ref{thm:F_0}), and \(Q = \{13,23,24,15\}\). In particular, among the
    minimal obstructions in the characterization in Theorem~\ref{thm:charac}, \(Q\) has
    the weakest known superlinear lower bound. Improving this bound would therefore
    improve the universal lower bound for ordered trees with superlinear Ramsey numbers.
    \item The following construction also gives a superlinear lower bound on
    \(R_<(F,K_n)\) for \(F = F_2 = \{12,23,14\}\) and \(F=S:=\{14,23,26,45\}\). Set \(L
    := \lceil \log_2(n)\rceil\) and for each vertex \(x\in[n]\), choose \(\lambda(x) \in
    \{0,1,\ldots,L- \nolinebreak 1\}\) uniformly at random and independently. Define an
    ordered graph \(G\) on \([n]\) as follows: for \(x<y\), let \(xy\in E(G)\) if and
    only if \(\lambda(x)\le\lambda(y) \quad\text{and}\quad 2^{\lambda(x)}\le
    y-x<2^{\lambda(x)+1}\). This graph is \(F_2\)-free and \(S\)-free and one can show
    that w.h.p.~\(\alpha(G) = O(\frac{n\log\log n}{\log n})\). This does slightly worse
    than Construction \ref{const:F1,F2} due to the \(\log\log n\) term in the numerator
    (cf. Theorem~\ref{thm:F1,F2 independence number}).

    \item By careful utilization of the stepping-up lemma of Erd\H{o}s and Hajnal,
    constructions for the off-diagonal ordered Ramsey numbers for finite
    families~\(R_<(\mathcal{F},K_n)\) may be used to obtain lower bounds for the
    (unordered) off-diagonal hypergraph Ramsey numbers~\(R(H,K_n^{(3)})\). For more
    information on the implicit relationship between the family~\(\mathcal{F}\) of
    ordered graphs and the~\(3\)-graph~\(H\), see \cite{CFGHMSV}, where this method was
    developed, or \cite{HNWY}, where it was lifted to uniformities larger than~\(3\). We
    wonder if our constructions here can be used, via this framework, to obtain new
    lower bounds on off-diagonal hypergraph Ramsey numbers.
\end{itemize}

\paragraph{Acknowledgment.} We thank Dingding Dong and Dylan King for many useful
discussions on this topic, Xiaoyu He for telling us about the connection to hypergraph
Ramsey numbers, and Ruben Ascoli for reviewing an earlier version of this manuscript.

\paragraph{AI Declaration.} Most constructions and proofs in this paper were obtained by
the authors without the use of AI. ChatGPT-6 Astra was used substantially in the
following parts: the structural arguments in Section~\ref{sec:charac}; the proofs in
Section~\ref{sec:Ramsey-good}; the construction used for the lower bound in
Theorem~\ref{thm:F4} (the authors previously had a construction similar to
Construction~\ref{const:randomized-Tutte}, yielding a bound with an extra \(\log\log n\)
factor in the denominator); the analysis of Construction~\ref{const:DS} (the
construction itself was devised by the authors, but the proof of
Lemma~\ref{lem:DS construction independence number} was inspired by a proof generated by
ChatGPT for a different construction); the proofs of
Theorems~\ref{thm:one-vertex interlace} and~\ref{thm:chi-bounding reduction} (the
proof ideas are entirely due to the authors, who asked ChatGPT to generalize their
arguments for the representative cases \(\{12,14,35\}\) and \(\{14,23,35\}\)); and a
simplification of the analysis of Construction~\ref{const:randomized-Tutte}. The authors
take full responsibility for the correctness of all proofs.

\fontsize{11pt}{12pt}
\selectfont
\hypersetup{linkcolor={red!70!black}}
\setlength{\parskip}{2pt plus 0.3ex minus 0.3ex}

\bibliographystyle{auxfile.bst}
\bibliography{bib.bib}

@article{FurediHajnal,
  author  = {F{\"u}redi, Zolt{\'a}n and Hajnal, P{\'e}ter},
  title   = {{D}avenport--{S}chinzel theory of matrices},
  journal = {Discrete Mathematics},
  volume  = {103},
  number  = {3},
  pages   = {233--251},
  year    = {1992},
  doi     = {10.1016/0012-365X(92)90316-8}
}

@article{PachTardos,
  author  = {Pach, J{\'a}nos and Tardos, G{\'a}bor},
  title   = {Forbidden paths and cycles in ordered graphs and matrices},
  journal = {Israel Journal of Mathematics},
  volume  = {155},
  number  = {1},
  pages   = {359--380},
  year    = {2006},
  doi     = {10.1007/BF02773960}
}

@article{CFLS,
  author  = {Conlon, David and Fox, Jacob and Lee, Choongbum and Sudakov, Benny},
  title   = {Ordered {R}amsey numbers},
  journal = {Journal of Combinatorial Theory, Series B},
  volume  = {122},
  pages   = {353--383},
  year    = {2017},
  doi     = {10.1016/j.jctb.2016.06.007}
}

@article{BCKK,
  author  = {Balko, Martin and Cibulka, Josef and Kr{\'a}l, Karel and Kyn{\v{c}}l, Jan},
  title   = {{R}amsey numbers of ordered graphs},
  journal = {The Electronic Journal of Combinatorics},
  volume  = {27},
  number  = {1},
  pages   = {P1.16},
  year    = {2020},
  doi     = {10.37236/7816}
}

@incollection{Balko_survey,
  author    = {Balko, Martin},
  title     = {A survey on ordered {R}amsey numbers},
  booktitle = {Sum(m)it280: Surveys in Extremal Combinatorics and Combinatorial Geometry},
  editor    = {Katona, Gyula O. H. and Patk{\'o}s, Bal{\'a}zs and Tompkins, Casey},
  series    = {Bolyai Society Mathematical Studies},
  volume    = {32},
  pages     = {71--96},
  publisher = {Springer},
  address   = {Cham},
  year      = {2026},
  doi       = {10.1007/978-3-032-18810-6_4}
}

@article{MubayiSuk,
  author  = {Mubayi, Dhruv and Suk, Andrew},
  title   = {{R}amsey numbers of cliques versus monotone paths},
  journal = {European Journal of Combinatorics},
  volume  = {118},
  pages   = {103922},
  year    = {2024},
  doi     = {10.1016/j.ejc.2024.103922}
}

@article{GJS,
  author  = {Gishboliner, Lior and Jin, Zhihan and Sudakov, Benny},
  title   = {{R}amsey problems for monotone paths in graphs and hypergraphs},
  journal = {Combinatorica},
  volume  = {44},
  number  = {3},
  pages   = {509--529},
  year    = {2024},
  doi     = {10.1007/s00493-024-00082-7}
}

@article{GJJ,
  author  = {Gir{\~a}o, Ant{\'o}nio and Janzer, Barnab{\'a}s and Janzer, Oliver},
  title   = {Ordered {R}amsey numbers of powers of paths},
  journal = {Combinatorica},
  volume  = {46},
  number  = {4},
  pages   = {27},
  year    = {2026},
  doi     = {10.1007/s00493-026-00223-0}
}

@article{BJV,
  author  = {Balko, Martin and Jel{\'\i}nek, V{\'\i}t and Valtr, Pavel},
  title   = {On ordered {R}amsey numbers of bounded-degree graphs},
  journal = {Journal of Combinatorial Theory, Series B},
  volume  = {134},
  pages   = {179--202},
  year    = {2019},
  doi     = {10.1016/j.jctb.2018.06.002}
}

@article{ARU,
  author  = {Axenovich, Maria and Rollin, Jonathan and Ueckerdt, Torsten},
  title   = {Chromatic number of ordered graphs with forbidden ordered subgraphs},
  journal = {Combinatorica},
  volume  = {38},
  number  = {5},
  pages   = {1021--1043},
  year    = {2018},
  doi     = {10.1007/s00493-017-3593-0}
}

@inproceedings{Interval_minor,
  author    = {Bourneuf, Romain and Cocquet, Julien and Tang, Chaoliang and Thomass{\'e}, St{\'e}phan},
  title     = {A polynomial-time approximation algorithm for complete interval minors},
  booktitle = {Approximation, Randomization, and Combinatorial Optimization. Algorithms and Techniques (APPROX/RANDOM 2025)},
  editor    = {Ene, Alina and Chattopadhyay, Eshan},
  series    = {Leibniz International Proceedings in Informatics (LIPIcs)},
  volume    = {353},
  pages     = {15:1--15:23},
  publisher = {Schloss Dagstuhl -- Leibniz-Zentrum f{\"u}r Informatik},
  address   = {Dagstuhl, Germany},
  year      = {2025},
  doi       = {10.4230/LIPIcs.APPROX/RANDOM.2025.15}
}

@article{FKMV,
  author  = {F{\"u}redi, Zolt{\'a}n and Kostochka, Alexandr and Mubayi, Dhruv and Verstra{\"e}te, Jacques},
  title   = {Ordered and convex geometric trees with linear extremal function},
  journal = {Discrete \& Computational Geometry},
  volume  = {64},
  number  = {2},
  pages   = {324--338},
  year    = {2020},
  doi     = {10.1007/s00454-019-00149-z}
}

@inproceedings{PettieTardosWalczak,
  author    = {Pettie, Seth and Tardos, G{\'a}bor and Walczak, Bartosz},
  title     = {On a clique game and the {E}rd{\H{o}}s--{H}ajnal problem on high-chromatic high-girth subgraphs},
  booktitle = {Proceedings of the 2026 Annual ACM-SIAM Symposium on Discrete Algorithms (SODA)},
  editor    = {Larsen, Kasper Green and Saha, Barna},
  pages     = {2903--2927},
  publisher = {Society for Industrial and Applied Mathematics},
  address   = {Philadelphia, PA},
  year      = {2026},
  doi       = {10.1137/1.9781611978971.108}
}

@book{clrs,
  author    = {Cormen, Thomas H. and Leiserson, Charles E. and Rivest, Ronald L. and Stein, Clifford},
  title     = {Introduction to Algorithms},
  edition   = {4th},
  publisher = {MIT Press},
  year      = {2022},
  isbn      = {9780262046305}
}

@article{Alon,
  author  = {Alon, Noga},
  title   = {Independence numbers of locally sparse graphs and a {R}amsey type problem},
  journal = {Random Structures \& Algorithms},
  volume  = {9},
  number  = {3},
  pages   = {271--278},
  year    = {1996},
  doi     = {10.1002/(SICI)1098-2418(199610)9:3<271::AID-RSA1>3.0.CO;2-U}
}

@article{CFGHMSV,
  author  = {Conlon, David and Fox, Jacob and Gunby, Benjamin and He, Xiaoyu and Mubayi, Dhruv and Suk, Andrew and Verstra{\"e}te, Jacques},
  title   = {On off-diagonal hypergraph {R}amsey numbers},
  journal = {International Mathematics Research Notices},
  volume  = {2025},
  number  = {11},
  pages   = {rnaf122},
  year    = {2025},
  doi     = {10.1093/imrn/rnaf122}
}

@article{HNWY,
  author  = {He, Xiaoyu and Nie, Jiaxi and Wigderson, Yuval and Yu, Hung-Hsun},
  title   = {Off-diagonal {R}amsey numbers for linear hypergraphs},
  journal = {Combinatorics, Probability and Computing},
  volume  = {35},
  number  = {4},
  pages   = {554--567},
  year    = {2026},
  doi     = {10.1017/S0963548326100443}
}

@misc{GL,
  author        = {Gishboliner, Lior and Li, Xiangyu},
  title         = {Upper bounds for ordered {R}amsey numbers of forests and bounded-degree graphs},
  howpublished  = {arXiv preprint arXiv:2609.38828},
  year          = {2026},
  eprint        = {2609.38828},
  archivePrefix = {arXiv},
  primaryClass  = {math.CO},
  url           = {https://arxiv.org/abs/2609.38828}
}

@article{Chvatal,
  author  = {Chv{\'a}tal, Va{\v{s}}ek},
  title   = {Tree-complete graph {R}amsey numbers},
  journal = {Journal of Graph Theory},
  volume  = {1},
  number  = {1},
  pages   = {93},
  year    = {1977},
  doi     = {10.1002/jgt.3190010118}
}

@article{BE,
  author  = {Burr, Stefan A. and Erd{\H{o}}s, Paul},
  title   = {Generalizations of a {R}amsey-theoretic result of {C}hv{\'a}tal},
  journal = {Journal of Graph Theory},
  volume  = {7},
  number  = {1},
  pages   = {39--51},
  year    = {1983},
  doi     = {10.1002/jgt.3190070106}
}

@article{BP,
  author  = {Balko, Martin and Poljak, Marian},
  title   = {On off-diagonal ordered {R}amsey numbers of nested matchings},
  journal = {Discrete Mathematics},
  volume  = {346},
  number  = {2},
  pages   = {113223},
  year    = {2023},
  doi     = {10.1016/j.disc.2022.113223}
}

@unpublished{BDW,
  author = {Bria{\'n}ski, Marcin and Davies, James and Walczak, Bartosz},
  title  = {Coloring graphs excluding an induced ordered matching},
  note   = {In preparation}
}

@article{AAC,
  author  = {Aboulker, Pierre and Aubian, Guillaume and Charbit, Pierre},
  title   = {Heroes in oriented complete multipartite graphs},
  journal = {Journal of Graph Theory},
  volume  = {105},
  number  = {4},
  pages   = {652--669},
  year    = {2024},
  doi     = {10.1002/jgt.23061}
}

@article{HartSharir,
  author  = {Hart, Sergiu and Sharir, Micha},
  title   = {Nonlinearity of {D}avenport--{S}chinzel sequences and of generalized path compression schemes},
  journal = {Combinatorica},
  volume  = {6},
  number  = {2},
  pages   = {151--177},
  year    = {1986},
  doi     = {10.1007/BF02579170}
}

@article{Nivasch,
  author  = {Nivasch, Gabriel},
  title   = {Improved bounds and new techniques for {D}avenport--{S}chinzel sequences and their generalizations},
  journal = {Journal of the ACM},
  volume  = {57},
  number  = {3},
  pages   = {17:1--17:44},
  year    = {2010},
  doi     = {10.1145/1706591.1706597}
}

\end{document}